\documentclass{article}
\usepackage{amsmath,amsfonts,amsthm,amssymb,amscd,bm,bbm}
\usepackage[hidelinks]{hyperref}
\usepackage[usenames,dvipsnames]{color}
\usepackage{url}
\hypersetup{colorlinks, citecolor=ForestGreen, linkcolor=MidnightBlue, urlcolor=Blue}

\newcommand{\be}{\begin{equation}}
\newcommand{\ee}{\end{equation}}
\newcommand{\bp}{\begin{proof}}
\newcommand{\ep}{\end{proof}}
\newcommand{\bi}{\begin{itemize}}
\newcommand{\ei}{\end{itemize}}
\newcommand{\om}{\omega}
\newcommand{\omm}{\Omega}

\newcommand{\bbb}{{\cal B}}
\newcommand{\BBB}{\mathfrak B}
\newcommand{\ccc}{{\cal C}}
\newcommand{\dd}{{\cal D}}
\newcommand{\ees}{{\cal E}}
\newcommand{\fff}{{\cal F}}

\newcommand{\FFF}{{\bm {\mathcal{F}}}}
\newcommand{\hh}{\bm {\mathcal H}}

\newcommand{\RR}{\bm {\mathcal R}}
\newcommand{\SSS}{\bm {S}}
\newcommand{\gi}{{\cal G}}
\newcommand{\h}{{\cal H}}

\newcommand{\kkk}{{\cal K}}
\newcommand{\elll}{{\cal L}}
\newcommand{\mmm}{{\cal M}}
\newcommand{\nnn}{{\cal N}}

\newcommand{\ppp}{{\cal P}}
\newcommand{\PPPP}{\mathfrak P}
\newcommand{\qqq}{{\cal Q}}
\newcommand{\R}{{\cal R}}

\newcommand{\vvv}{{\cal V}}

\newcommand{\yy}{\mathbf{y}}

\newcommand{\e}{\mathbb{E}}

\newcommand{\GG}{\mathbb{G}}

\newcommand{\pp}{\mathbb{P}}

\newcommand{\rr}{\mathbb{R}}

\newcommand{\PP}{\mathsf{P}}
\newcommand{\QQ}{\mathsf{Q}}

\newcommand{\kp}{\varkappa}

\newcommand{\q}{\quad}

\newcommand{\f}{\frac}
\newcommand{\lm}{\lambda}
\newcommand{\p}{\partial}
\newcommand{\ph}{\varphi}
\newcommand{\bg}{\Bigr}

\newcommand{\De}{\delta}
\newcommand{\de}{\Delta}
\newcommand{\g}{\nabla}
\newcommand{\dt}{\dot}

\newcommand{\es}{\varepsilon}

\newcommand{\al}{\alpha}

\newcommand{\nt}{\noindent}

\newcommand{\ch} {\mathbbm{1}}
\newcommand{\tx} {\text}
 \newcommand{\iin} {\infty}
 
  \newcommand{\vo}{\varrho}
  \newcommand{\mpp}{\emph}
\theoremstyle{plain}
\newtheorem{theorem}{Theorem}[section]
\newtheorem*{mt}{Main Theorem}
\newtheorem{lemma}[theorem]{Lemma}
\newtheorem{proposition}[theorem]{Proposition}
\newtheorem{corollary}[theorem]{Corollary}
\theoremstyle{definition}
\newtheorem{definition}[theorem]{Definition}

\theoremstyle{remark}

\numberwithin{equation}{section}
\newtheorem*{definition*}{Definition}
\newtheorem*{problem*}{Problem}

\newtheorem*{remark*}{Remark}
\newtheorem*{note*}{Note}
\begin{document}
\author{Davit Martirosyan\footnote{Department of Mathematics, University of Cergy-Pontoise, CNRS UMR 8088, 2 avenue
Adolphe Chauvin, 95300 Cergy-Pontoise, France; e-mail: Davit.Martirosyan@u-cergy.fr}}
\title{Exponential mixing for the white\,-\,forced damped nonlinear wave equation} 

\maketitle

\begin{abstract}
The paper is devoted to studying the stochastic nonlinear wave (NLW) equation 
$$
\p_t^2 u+\gamma \p_t u-\de u+f(u)=h(x)+\eta(t,x)
$$
in a bounded domain $D\subset\rr^3$. The equation is supplemented with the Dirichlet boundary condition. Here $f$ is a nonlinear term, $h(x)$ is a function in $H^1_0(D)$ and $\eta(t,x)$ is a non-degenerate white noise. We show that the Markov process associated with the flow $\xi_u(t)=[u(t),\dt u(t)]$ has a unique stationary measure $\mu$, and the law of any solution converges to $\mu$ with exponential rate in the dual-Lipschitz norm.

\smallskip
\nt

\smallskip
\nt
\end{abstract}

\tableofcontents

\medskip

\nt 
\bigskip
\section{Introduction}
We consider the stochastic NLW equation
\be\label{1.1}
\p_t^2u+\gamma \p_tu-\de u+f(u)=h(x)+\eta(t,x),\q [u(0),\dt u(0)]=[u_0, u_1]
\ee
in a bounded domain $D\subset\rr^3$ with a smooth boundary. The equation is supplemented with the Dirichlet boundary condition. The nonlinear term $f$ satisfies the dissipativity and growth conditions that are given in the next section (see \eqref{1.8}-\eqref{1.6}). Here we only mention that they hold for functions  $f(u)=\sin u$ and $f(u)=|u|^{\rho}u-\lm u$, where $\lm$ and $\rho\in(0,2)$ are some constants. These functions correspond to the damped sine-Gordon and Klein-Gordon equations, respectively.
The force $\eta(t)$ is a white noise of the form
\be\label{1.9}
\eta(t,x)=\sum_{j=1}^\infty b_j\dt\beta_j(t)e_j(x).
\ee
Here $\{\beta_j(t)\}$ is a sequence of independent standard Brownian motions, $\{e_j\}$ is an orthonormal basis in $L^2(D)$ composed of the eigenfunctions of the Dirichlet Laplacian, and $\{b_j\}$ is a sequence of positive numbers that goes to zero sufficiently fast (see \eqref{2.24}). 
The initial point $[u_0,u_1]$ belongs to the phase space $\h=H^1_0(D)\times L^2(D)$. Finally, $h(x)$ is a function in $H^1_0(D)$. 
The following theorem is the main result of this paper.
\begin{mt}\label{1.4}
Under the above hypotheses, the Markov process associated with the flow $y(t)=[u(t),\dt u(t)]$ of equation \eqref{1.1} possesses a unique stationary measure $\mu\in\ppp(\h)$. Moreover, there are positive constants $C$ and $\kp$ such that
\be\label{6.53}
|\e \psi(y(t))-\int_{\h}\psi(z)\mu(dz)|\leq Ce^{-\kp t}\exp(\kp|y|_\h^4),\q t\geq 0, 
\ee
for any 1-Lipschitz function $\psi:\h\to\rr$, and any initial point $y\in\h$. 
\end{mt}

\nt
Thus, the limit of the average of $\psi(y(t))$ is a quantity that does not depend on the initial point.

Before outlining the main ideas of the proof of this result, let us discuss some of the earlier works concerning the ergodicity of the stochastic nonlinear PDE's and the main difficulties that occur in our case. In the context of stochastic PDE's, the initial value problem and existence of a stationary measure was studied by Vishik--Fursikov--Komech \cite{VKF-1979} for the stochastic Navier--Stokes system and later developed for many other problems (see the references in \cite{DZ1992}). The uniqueness of stationary measure and its ergodicity are much more delicate questions. First results in this direction were obtained in the papers \cite{FM-1995, KS-cmp2000, EMS-2001, BKL-2002} devoted to the Navier--Stokes system and other PDE's arising in mathematical physics (see also \cite{MR1245306, MR1641664} and Part III in \cite{DZ1996} for some 1D parabolic equations). They were later extended to equations with multiplicative and very degenerate noises \cite{odasso-2008,HM-2008}. We refer the reader to the recent book \cite{KS-book} and the review paper \cite{debussche2013ergodicity} for a detailed account of the main results obtained so far. 

We now discuss in more details the case of dispersive equations, for which fewer results are known.
One of the first results on the ergodicity of dispersive PDE's was stablished in the paper of E, Khanin, Mazel and Sinai \cite{EWKMS2000}, where the authors prove the existence and uniqueness of  stationary measure for the one dimensional inviscid Burgers equation perturbed by a space-periodic white noise. The qualitative study of stationary solutions is also carried out, and the analysis relies on the Lax-Oleinik variational principle.
The ergodicity of a white-forced NLW equation was studied by Barbu and Da Prato \cite{BD-2002}, where the authors prove the existence of stationary distribution for a nonlinearity which is a non-decreasing function satisfying the growth restriction $|f''(u)|\leq C(|u|+1)$, and some standard dissipativity conditions. Uniqueness is established under the additional hypotheses, that $f$ satisfies \eqref{1.8} with $\rho<2$, and $\sup\{|f'(u)|\cdot |u|^{-\rho}, u\in\rr\}$ is sufficiently small. In the paper by Debussche and Odasso \cite{DO-2005}, the authors establish the convergence to the equilibrium with polynomial speed at any order (\mpp{polynomial mixing}) for weakly damped nonlinear Schr\"odinger equation. The proof of this result relies on the coupling argument. The main difficulty in establishing the exponential rate of convergence is due to the complicated Lyapunov structure and the fact that the Foa\c{s}-Prodi estimates hold in average and not path-wise. 
In \cite{DirSoug2005}, Dirr and Souganidis study the Hamilton-Jacobi equations perturbed by additive noise. They show, in particular, that under suitable assumptions on the Hamiltonian, the stochastic equation has a unique up to constants space-periodic global attracting solution, provided the unperturbed equation possesses such solution.  
In the recent paper by Debussche and Vovelle \cite{debussche2013invariant} the existence and uniqueness of stationary measure is studied for scalar periodic first-order conservation laws with additive noise in any space dimension. It generalizes to higher dimensions the results established in \cite{EWKMS2000} (see also \cite{IturK2003}). In another recent paper \cite{BCK2014} by Bakhtin, Cator and Khanin, the authors study the ergodicity of the Burgers equation perturbed by a space-time stationary random force. It is proved, in particular, that the equation possesses space-time stationary global solutions, and that they attract all other solutions. The proof uses the Aubry-Mather theory for action-minimizing trajectories, and weak KAM theory for the Hamilton-Jacobi equations. 

In the present paper we extend the results established in \cite{BD-2002}, proving that the hypotheses  $f'\geq 0$ and $\sup\{|f'(u)|\cdot |u|^{-\rho}, u\in\rr\}$ is small are not needed, and that the convergence to the equilibrium has exponential rate. We also show that the conclusion of the Main Theorem remains true for a force that is non-degenerate only in the low Fourier modes (see Theorem \ref{6.10}). The proof mainly relies on the coupling argument.
 
Of course, one of the main difficulties when dealing with dispersive PDE's comes from the lack of the regularizing property, and with it, of some well-known compactness arguments. As a consequence, this changes the approach when showing the stability of solutions. In particular, this is the case, when establishing the Foia\c{s}-Prodi estimate for NLW (Proposition \ref{4.13}). Moreover, this estimate (which shows that the large time behavior of solutions is determined by finitely many modes and enables one to use the Girsanov theorem) differs from the classical one, since the growth of the intermediate process should be controlled (see inequality \eqref{4.16}). Due to the last fact, the  coupling constructed through the projections of solutions (cf. \cite{shirikyan-bf2008, odasso-2008}) does not ensure exponential rate of convergence. We therefore introduce a new type of coupling constructed via the intermediate process (see \eqref{2.25}-\eqref{6. Coupling relation 1}).
The same difficulty occurs when showing the recurrence of solutions, i.e. that the trajectory of the solution enters arbitrarily small ball with positive probability in a finite time (Proposition \ref{6.12}). The standard argument to show this property is the use of the portmanteau theorem. However, due to the lack of the smoothing effect, the portmanteau technique is not applicable, and another approach is proposed. 

Without going into details, we give an informal description of our approach. The proof of the existence of stationary measure is rather standard and relies on the Bogolyubov-Krylov argument, which ensures the existence, provided the process $y(t)=[u(t),\dt u(t)]$ has a uniformly bounded moment in some $\h$-compact space. To obtain such a bound, we follow a well-known argument coming from the theory of attractors (e.g., see \cite{BV1992, Har85}). Namely, we split the function $u$ to the sum $u=v+z$, where, roughly speaking, $v$ takes the Brownian of equation, and $z$-nonlinearity. We then show that the corresponding flows have uniformly bounded moments in $\h^s=H^{1+s}(D)\times H^s(D)$ for $s>0$ sufficiently small (Proposition \ref{2.11}). The bound for $|[v(t),\dt v(t)]|_{\h^s}$ follows from the It\^o formula, while that of $|[z(t),\dt z(t)]|_{\h^s}$ is based on the argument similar to the one used in \cite{zelik2004}.
The proof of exponential mixing relies on Theorem 3.1.7 in \cite{KS-book}, which gives a general criterion that ensures the convergence to the equilibrium with exponential rate. Construction of a coupling that satisfies the hypotheses of the mentioned theorem is based on four key ingredients: the Foia\c{s}-Prodi estimate for NLW, the Girsanov theorem, the recurrence property of solutions, and the stopping time technique. 

Finally, we make some comments on the hypotheses imposed on the nonlinear term $f$ and the coefficients $b_j$ entering the definition of the force $\eta$. Inequalities \eqref{1.5}-\eqref{1.6} are standard in the study of NLW equation, they ensure that the Cauchy problem is well-posed (e.g., see \cite{CV2002} and \cite{Lions1969} for deterministic cases). The hypothesis $\rho<2$ is needed to prove the stability of solutions. The fact that the coefficients $b_j$ are not zero ensures that $\eta$ is non-degenerate in all Fourier modes, which is used to establish the recurrence of solutions and exponential squeezing. As was mentioned above, we show that this condition could be relaxed.

The paper is organized as follows. In Section \ref{2.26} we announce the main result and outline the scheme of its proof. Next, the large time behavior and stability of solutions are studied in Sections \ref{3.0} and \ref{4.0}, respectively. Finally, the complete proof of the main result is presented in Section \ref{6.31}.

\bigskip
{\bf Acknowledgments}. I am grateful to my supervisor Armen Shirikyan, for attracting my attention to this problem, and for many fruitful discussions. This research was carried out within the MME-DII Center of Excellence (ANR 11 LABX 0023 01) and partially supported by the ANR grant STOSYMAP (ANR 2011 BS01 015 01).

\subsection*{Notation}
For an open set $D$ of a Euclidean space and separable Banach spaces $X$ and $Y$, we introduce the following function spaces:

\nt
$L^p=L^p(D)$ is the Lebesgue space of measurable functions whose $p^{\text{th}}$ power is integrable. In the case $p=2$ the corresponding norm is denoted by $\|\cdot\|$.

\nt
$H^s=H^s(D)$ is the Sobolev space of order $s$ with the usual norm $\|\cdot\|_s$.

\nt
$H^s_0=H^s_0(D)$ is the closure in $H^s$ of infinitely smooth functions with compact support.

\nt
$H^{1,p}=H^{1,p}(D)$ is the Sobolev space of order $1$ with exponent $p$, that is, the space of $L^p$ functions whose first order derivatives remain in $L^p$.

\nt
$L(X,Y)$ stands for the space of linear continuous operators from $X$ to $Y$ endowed with the natural norm.

\nt
$C_b(X)$ is the space of continuous bounded functions $\psi:X\to\rr$ endowed with the norm of uniform convergence:
$$
|\psi|_\infty=\sup_{x\in X}|\psi(x)|.
$$
$L_b(X)$ is the space of bounded Lipschitz functions, i.e. of functions $\psi\in C_b(X)$ such that
$$
|\psi|_L:=|\psi|_\infty+\sup_{x\neq y}\f{|\psi(x)-\psi(y)|}{|x-y|_X}<\infty.
$$

\nt
$B_X(R)$ stands for the ball in $X$ of radius $R$ and centered at the origin.

\nt
$\bbb(X)$ is the Borel $\sigma$-algebra of subsets of $X$.

\nt
$\ppp(X)$ denotes the space of probability Borel measures on $X$. Two metrics are defined on the space $\ppp(X)$: the metric of total variation
$$
|\mu_1-\mu_2|_{var}=\sup_{\Gamma\in \bbb(X)}|\mu_1(\Gamma)-\mu_2(\Gamma)|,
$$
and the dual Lipschitz metric
$$
|\mu_1-\mu_2|_L^*=\sup_{|\psi|_L\leq 1}|(f,\mu_1)-(f,\mu_2)|,
$$
where $(\psi,\mu)$ denotes the integral of $\psi$ over $X$ with respect to $\mu$. 

\nt
Finally, by $C_1,C_2,\ldots$, we shall denote unessential positive constants.

\section{Exponential mixing}\label{2.26}
We start this section by a short discussion of the well-posedness of the Cauchy problem for equation \eqref{1.1}. We then state the main result and outline the scheme of its proof.
\subsection{Existence and uniqueness of solutions}
Before giving the definition of a solution of equation \eqref{1.1}, let us make the precise hypotheses on the nonlinearity and the coefficients entering the definition of $\eta(t)$. We suppose that the function $f$ satisfies the growth restriction
\be\label{1.8}
|f''(u)|\leq C(|u|^{\rho-1}+1),\q u\in\rr,
\ee  
where $C$ and $\rho<2$ are positive constants, and the dissipativity conditions
\begin{align}
F(u)&\geq -\nu u^2-C,\q u\in\rr\label{1.5},\\
f(u)u- F(u)&\geq-\nu u^2-C, \q u\in\rr\label{1.6},
\end{align}
where $F$ is the primitive of $f$, $\nu\leq (\lm_1\wedge\gamma)/8$ is a positive constant, and $\lm_j$ stands for the eigenvalue corresponding to $e_j$. The coefficients $b_j$ are supposed to be positive numbers satisfying
\be \label{2.24}
\BBB=\sum_{j=1}^\iin b_j^2<\iin,\q\BBB_1=\sum_{j=1}^\iin\lm_j b_j^2<\infty.
\ee
Let us introduce the functions
$$
g_j=[0,b_j e_j],\q \hat\zeta(t)=\sum_{j=1}^\iin\beta_j(t)g_j.
$$
\begin{definition}\label{definition 2.1}
Let $y_0=[u_0,u_1]$ be a $\h$-valued random variable defined on a complete probability space $(\omm,\fff,\pp)$ that is independent of $\hat\zeta(t)$. A random process $y(t)=[u(t),\dt u(t)]$ defined on  $(\omm,\fff,\pp)$ is called \emph{a solution} (or \emph{a flow}) of equation \eqref{1.1} if the following two  conditions hold:
\bi 
\item Almost every trajectory of $y(t)$ belongs to the space $C(\rr_+;\h)$, and the process $y(t)$ is adapted to the filtration $\fff_t$ generated by $y_0$ and $\hat\zeta(t)$.

\item Equation \eqref{1.1} is satisfied in the sense that, with probability 1,
\be\label{2.1}
y(t)=y_0+\int_0^t g(s)\,ds+\hat\zeta(t),\q t\geq 0,
\ee
where we set
$$
g(t)=[\dt u,-\gamma \dt u+\de u-f(u)+h(x)],
$$
and relation \eqref{2.1} holds in $L^2\times H^{-1}$.
\ei
\end{definition}

\nt
Let us endow the space $\h$ with the norm
$$
|y|_{\h}^2=\|\g y_1\|^2+\|y_2+\al y_1\|^2\q \text{ for } y=[y_1,y_2]\in\h,
$$
where $\al>0$ is a small parameter.
Introduce the energy functional
\be 
\ees(y)=|y|_\h^2+2\int_D F(y_1)\,dx, \q y=[y_1,y_2]\in\h,
\ee 
and let $\ees_u(t)=\ees(y(t))$. We have the following theorem.
\begin{theorem}\label{2.12}
Under the above hypotheses, let $y_0$ be an $\h-$valued random variable that is independent of $\hat\zeta$ and satisfies $\e\ees(y_0)<\iin$. Then equation \eqref{1.1} possesses a solution in the sense of Definition \ref{definition 2.1}. Moreover, it is unique, in the sense that if $\tilde y(t)$ is another solution, then with $\pp$-probability 1 we have $y(t)=\tilde y(t)$ for all $t\geq 0$. In addition, we have the a priori estimate
\be\label{2.2}
\e\ees_u(t)\leq\e\ees_u(0)e^{-\al t}+C(\gamma,\BBB,\|h\|).
\ee
\end{theorem}
We refer the reader to the book \cite{DZ1992} for proofs of similar results. We confine ourselves to the formal derivation of inequality \eqref{2.2} in the next section.
\subsection{Main result and scheme of its proof}\label{Main result and scheme of its proof}
Let us denote by $S_t(y,\cdot)$ the flow of equation \eqref{1.1} issued from the initial point $y\in\h$. A standard argument shows that $S_t(y,\cdot)$ defines a Markov process in $\h$ (e.g., see \cite{DZ1992, KS-book}). We shall denote by $(y(t),\pp_y)$ the corresponding Markov family. In this case, the Markov operators have the form

\begin{align*}
\PPPP_t\psi(y)&=\int_\h\psi(z)P_t(y,dz) \q \text{ for any } \psi\in C_b(\h),\\
\PPPP^*_t\lm(\Gamma)&=\int_\h P_t(y,\Gamma)\lm(dy)\q \text{ for any } \lm\in\ppp(\h),
\end{align*}
where $P_t(y,\Gamma)=\pp_y(S_t(y,\cdot)\in\Gamma)$ is the transition function.
The following theorem on exponential mixing is the main result of this paper.
\begin{theorem}\label{2.15}
Under the hypotheses of Theorem \ref{2.12}, the Markov process associated with the flow of equation \eqref{1.1} has a unique stationary measure $\mu\in\ppp(\h)$. Moreover, there exist positive constants $C$ and $\kp$ such that for any $\lm\in\ppp(\h)$ we have
\be\label{6. In-main inequality}
|\PPPP^*_t\lm-\mu|_L^*\leq C e^{-\kp t}\int_\h \exp(\kp|y|_\h^4)\,\lm(dy).
\ee
\end{theorem}

{\it Scheme of the proof.}
We shall construct an extension for the family $(y(t),\pp_y)$ that satisfies the hypotheses of Theorem 3.1.7 in \cite{KS-book}, providing a general criterion for exponential mixing. To this end, let us fix an initial point $\yy=(y,y')$ in $\hh=\h\times\h$, and let $\xi_u=[u,\p_t u]$ and $\xi_{u'}=[u',\p_t u']$ be the flows of equation \eqref{1.1} that are issued from $y$ and $y'$, respectively. Consider an intermediate process $v$, which is the solution of 
\be\label{2.25}
\p_t^2v+\gamma \p_t v-\de v+f(v)+P_N[f(u)-f(v)]=h(x)+\eta(t,x), \q\xi_v(0)=y'.
\ee
Let us denote by $\lm(y,y')$ and $\lm'(y,y')$ the laws of the processes $\{\xi_v\}_T$ and $\{\xi_{u'}\}_T$, respectively, where $\{z\}_T$ stands for the restriction of $\{z(t); t\geq 0\}$ to $[0,T]$. Thus, $\lm$ and $\lm'$ are probability measures on $C(0,T;\h)$. Let $(\vvv(y,y'),\vvv'(y,y'))$ be a maximal coupling for $(\lm(y,y'),\lm'(y,y'))$. By Proposition 1.2.28 in \cite{KS-book}, such a pair exists and can be chosen to be a measurable function of its arguments. For any $s\in [0,T]$, we shall denote by $\vvv_s$ and $\vvv'_s$ the restrictions of $\vvv$ and $\vvv'$ to the time $s$. Denote by $[\tilde v,\p_t \tilde v]$ and $[\tilde u',\p_t \tilde u']$ the corresponding flows. Then we have
\be\label{7.1}
\p_t^2 \tilde v+\gamma \p_t \tilde v-\de\tilde v+f(\tilde v)-P_N f(\tilde v)=h(x)+\psi(t),\q \xi_{\tilde v}(0)=y',
\ee
where $\psi$ satisfies
\be\label{6.30}
\dd\{\int_0^t \psi(s)\,ds\}_T=\dd\{\zeta(t)-\int_0^t P_N f(u)\,ds\}_T.
\ee
Introduce an auxiliary process $\tilde u$, which is the solution of 
\be 
\p_t^2 \tilde u+\gamma \p_t \tilde u-\de\tilde u+f(\tilde u)-P_N f(\tilde u)=h(x)+\psi(t),\q \xi_{\tilde u}(0)=y.
\ee
Let us note that $u$ satisfies the same equation, where $\psi$ should be replaced by $\eta(t)-P_N f(u)$. In view of \eqref{6.30}, we have (see the appendix for the proof)
\begin{equation}\label{6.40}
\dd\{\xi_{\tilde u}\}_T=\dd\{\xi_{u}\}_T.
\end{equation}
Introduce 
\be
\R_t(y,y')=\xi_{\tilde u}(t), \q \R_t'(y,y')=\xi_{\tilde u'}(t) \q \text{ for } t\in[0,T].\label{6. Coupling relation 1}
\ee
It is clear that $\RR_t=(\R_t,\R'_t)$ is an extension of $S_t(y)$ on the interval $[0,T]$. Let $\SSS_t=(S_t(\yy), S'_t(\yy))$ be the extension of $S_t(y)$ constructed by iteration of $\RR_t=(\R_t,\R'_t)$ on the half-line $t\geq 0$ (we do not recall here the procedure of construction, see the paper \cite{shirikyan-bf2008} for the details). With a slight abuse of notation, we shall keep writing $[\tilde u, \p_t\tilde u]$ and $[\tilde u', \p_t\tilde u']$ for the extensions of these two processes, and write $\xi_{\tilde v}(t)=\vvv_s(\SSS_{kT}(\yy))$ for $t=s+kT,\q 0\leq s<T$. This will not lead to a confusion.

\nt For any continuous process $y(t)$ with range in $\h$, we introduce the functional
\be\label{2.13}
\fff_y(t)=|\ees(y(t))|+\al\int_0^t|\ees (y(s))|\,ds,
\ee
and the stopping time
\be\label{6 In-6,1,1}
\tau_y=\inf\{t\geq 0:\fff_y(t)\geq \fff_y(0)+(L+M)t+r\},
\ee
where $L,M$ and $r$ are some positive constants to be chosen later.
In the case when $y$ is a process of the form $y=[z,\dt z]$, we shall write, $\fff^z$ and $\tau^z$ instead of $\fff_{[z,\dt z]}$ and $\tau_{[z,\dt z]}$, respectively.
Introduce the stopping times:
\begin{align*}
 \vo&=\inf\{t=s+kT: \vvv_s(\SSS_{kT}(\yy))\neq \vvv'_s(\SSS_{kT}(\yy))\}\equiv\inf\{t\geq 0: \xi_{\tilde v}(t)\neq \xi_{\tilde u'}(t)\},\\
 \tau&=\tau^{\tilde u}\wedge\tau^{\tilde u'},\q \sigma=\vo\wedge\tau.
 \end{align*}
Suppose that we are able to prove the following.
\begin{theorem}\label{2.14}
Under the hypotheses of Theorem \ref{2.15}, there are positive constants $\al,\De,\kp,d$ and $C$ such that the following properties hold.\\
\text{\textnormal{(Recurrence):}} For any $\yy=(y,y')\in\hh$, we have
\begin{align}
\e_y\exp(\kp\ees(y(t))&\leq\e_y\exp(\kp\ees (y(0))e^{-\al t}+C(\gamma,\BBB,\|h\|), \label{2.17}\\
\e_{\yy}\exp(\kp\tau_d)&\leq C(1+|\yy|_{\hh}^4), \label{2.16}
\end{align}
where $\tau_d$ stands for the first hitting time of the set $B_{\hh}(d)$.\\
\text{\textnormal{(Exponential squeezing):}} For any $\yy\in B_{\hh}(d)$, we have
\begin{align}
|S_t(\yy)-S_t'(\yy)|^2_\h&\leq Ce^{-\al t}|y-y'|_\h^2 \q\text{ for } 0\leq t\leq\sigma\label{6.33},\\
\pp_\yy\{\sigma=\infty\}&\geq\De\label{6.3},\\
\e_\yy[\ch_{\{\sigma<\infty\}}\exp(\De\sigma)]&\leq C\label{6.4},\\
\e_\yy[\ch_{\{\sigma<\infty\}} |\yy(\sigma)|_{\hh}^{8}]&\leq C\label{6.5}.
\end{align}
\end{theorem}
In view of Theorem 3.1.7 in \cite{KS-book}, this will imply Theorem \ref{2.15}. We establish Theorem \ref{2.14} in Section \ref{6.31}. The proof of \mpp{recurrence} relies on the Lyapunov function technique, while the proof of \mpp{exponential squeezing} is based on the Foia\c{s}-Prodi type estimate for equation \eqref{1.1}, the Girsanov theorem and the stopping time argument.

\subsection{Law of large numbers and central limit theorem}\label{Law of large numbers and central limit theorem}
Theorem \ref{2.15} implies the following result, which follows from inequality \eqref{6. In-main inequality} and some results established in Section 2 of \cite{shirikyan-ptrf2006}.
\begin{theorem}
Under the hypotheses of Theorem \ref{2.15}, for any Lipschitz bounded functional $\psi:\h\to\rr$ and any solution $y(t)=[u(t),\dt u(t)]$ of equation \eqref{1.1} issued from a non-random point $y_0\in\h$, the following statements hold.\\
\text{\textnormal{Strong law of large numbers.}} For any $\es>0$ there is an almost surely finite random constant $l\geq 1$ such that
\be\label{6.51}
|t^{-1}\int_0^t \psi(y(s))\,ds-(\psi,\mu)|\leq C(y_0,\psi) t^{-\f{1}{2}+\es} \q \text{ for } t\geq l.
\ee
\text{\textnormal{Central limit theorem.}} If $(\psi,\mu)=0$, there is a constant $a\geq 0$ depending only on $\psi$, such that for any $\es>0$, we have
\be\label{6.52}
\sup_{z\in\rr}(\theta_a(z)\cdot |\pp\{t^{-\f{1}{2}}\int_0^t\psi(y(s))\,ds\leq z\}-\Phi_a(z)|)\leq C(y_0,\psi)t^{-\f{1}{4}+\es},
\ee
where we set
$$
\theta_a(z)\equiv 1,\q \Phi_a(z)=\f{1}{a\sqrt{2\pi}}\int_{-\iin}^z e^{-\f{s^2}{2 a^2}}\,ds \q\text{ for } a>0,
$$
and 
$$
\theta_0(z)=1\wedge|z|,\q \Phi_0(z)=\ch_{\rr_+}(z).
$$
\end{theorem}

\nt
The proof of inequalities \eqref{6.51} and  \eqref{6.52} follow, respectively, from Corollary 3.4 and Theorem 2.8 in \cite{shirikyan-ptrf2006}, combined with inequalities \eqref{6.53} and \eqref{2.17}. \section{Large time estimates of solutions}\label{3.0}
The goal of this section is to analyze the dynamics of solutions and to obtain some a priori estimates for them.
\subsection{Proof of inequality \eqref{2.2}}

Let us apply the It\^o formula to the function $\GG(y)=|y|_\h^2$. Recall that for the process of the form \eqref{2.1}, the It\^o formula gives 
\be\label{2.3}
\GG(y(t))=\GG(y(0))+\int_0^t A(s)\,ds+\sum_{j=1}^\iin\int_0^t B_j(s)d\beta_j(s),
\ee
where we set
$$
A(t)=(\p_y \GG)(y(t);g(t))+\f{1}{2}\sum_{j=1}^\iin(\p_y^2 \GG)(y(t);g_j,g_j),\q
B_j(t)=(\p_y \GG)(y(t);g_j).
$$
Here $(\p_y\GG)(y;v)$ and $(\p_y^2\GG)(y;v,v)$ stand for the values of the first- and second-order derivatives of $\GG$ on the vector $v$.  
Since for $\GG(y)=|y|_\h^2$ we have
$$
\p_y\GG(y;\bar y)=2(y,\bar y)_\h,\q \p_y^2\GG(y;\bar y,\bar y)=2|\bar y|_\h^2,
$$
relation \eqref{2.3} takes the form
\be\label{2.3,2}
|y(t)|_{\h}^2=|y(0)|_\h^2+2\int_0^t(y,g)_\h\,ds+t\cdot\sum_{j=1}^\infty|g_j|_\h^2 +2\sum_{j=1}^\infty\int_0^t(y,g_j)_\h\,d\beta_j(s).
\ee
Let us note that
\begin{align}
(y,g)_\h&=(\g u,\g \dt u)+(\dt u+\al u,-\gamma \dt u+\de u-f(u)+h(x)+\al \dt u)\notag\\
&=-\al\|\g u\|^2-(\gamma-\al)\| \dt u\|^2+(\al^2-\al\gamma)(u,\dt u)+(\dt u+\al u,h)\notag\\
&\q\,-(\dt u+\al u,f(u))\label{2.4}.
\end{align}
By the Young and Poincar\'e inequalities, we have
\begin{align}
|(\al^2-\al\gamma)(u,\dt u)|&\leq \f{\al}{16}\|\g u\|^2+\f{4\al(\gamma-\al)^2}{\lm_1}\|\dt u\|^2\label{6.65},\\
|(\al u,h)|&\leq \f{\al}{16}\|\g u\|^2+\f{4\al}{\lm_1}\|h\|^2,\\
|(\dt u,h)|&\leq \f{\gamma-\al}{4}\|\dt u\|^2+(\gamma-\al)^{-1}\|h\|^2.\label{6.66}
\end{align}
Note also that, thanks to inequality \eqref{1.6}, we have 
\be\label{2.3,3}
-\al f(u)u\leq - \al F(u)+\al \nu u^2+\al C\leq - \al F(u)+\f{\al \lm_1}{8} u^2+\al C,
\ee
so that
\be\label{6.67}
-(\al u,f(u))\leq  - \al\int_D F(u)+\f{\al}{8} \|\g u\|^2+\al C\cdot\text{Vol}(D)
\ee
Now, by substituting \eqref{2.4} into \eqref{2.3,2}, using inequalities \eqref{6.65}-\eqref{6.67}, and noting that
$$
\int_0^t(\dt u,f(u))\,ds=\int_0^t\f{d}{ds}F(u(s))\,ds=F(u(t))-F(u(0)),
$$
we obtain that for $\al>0$ sufficiently small
\be\label{2.19}
\ees_u(t)\leq \ees_u(0)+\int_0^t\left(-\al\ees_u(s)+\kkk\right)\,ds-\f{\al}{2}\int_0^t|y(s)|_\h^2\,ds+M(t),
\ee
where $\kkk>0$ depends only on $\gamma, \BBB$ and $\|h\|$, and $M(t)$ is the stochastic integral
\be\label{2.20}
M(t)=2\sum_{j=1}^\infty b_j\int_0^t (\dt u+\al u,e_j)\,d\beta_j(s).
\ee
Taking the mean value in inequality \eqref{2.19} and using the Gronwall comparison principle, we arrive at \eqref{2.2}.
\subsection{Exponential moment of the flow}
In the following proposition we establish the uniform boundedness of exponential moment of $|\xi_u(t)|_\h$.
\begin{proposition}\label{proposition_exp}
Under the hypotheses of Theorem \ref{2.12}, there exists $\kp>0$ such that if the random variable $\ees_u(0)$ satisfies
$$
\e\exp(\kp\ees_u(0))<\infty,
$$
then
\be\label{2. In-exponential-moment-of-the-flow} 
\e\exp(\kp\ees_u(t))\leq\e\exp(\kp\ees_u(0))e^{-\al t}+C(\gamma,\BBB,\|h\|).
\ee
\end{proposition}
\bp
We represent $\xi_u(t)$ in the form \eqref{2.1}, and apply the It\^o formula \eqref{2.3} to the function  
$$
\GG(y)=\exp(\kp\ees(y)).
$$
Since 
\begin{align*}
\p_y\GG(y,\bar y)&=2\kp\,\GG(y)((y,\bar y)_\h+(f(y_1),\bar y_1)),\\
\p_y^2\GG(y;\bar y,\bar y)&=2\kp\,\GG(y)(2\kp\,((y,\bar y)_\h+(f(y_1),\bar y_1)^2+|\bar y|_\h^2+(f'(y_1),\bar y_1^2),
\end{align*}
we have
\begin{align*}
\p_y\GG(y;g)&=2\kp\,\GG(y)((y,g)_\h+(f(u),\dt u)),\\
\p_y^2\GG(y;g_j,g_j)&=2\kp\,\GG(y)(2\kp(y,g_j)_\h^2+|g_j|_\h^2).
\end{align*}
Hence, relation \eqref{2.3}, after taking the mean value, takes the form
$$
\e\GG(y(t))=\e\GG(y(0))+\kp\,\e\int_0^t \GG(y(s))\mmm(s)\,ds,
$$
where
$$
\mmm(t)=2((y,g)_\h+(f(u),\dt u))+2\kp\sum_{j=1}^\iin(y,g_j)_\h^2+\sum_{j=1}^\iin|g_j|_\h^2.
$$
Now note that by developing the expression $(y,g)_\h+(f(u),\dt u)$, the term $(f(u),\dt u)$ will disappear (see \eqref{2.4}). There remains another term containing $f$, namely the term $(-\al u,f(u))$, but this can be estimated using inequality \eqref{6.67}. Let us choose $\kp>0$ so small that $\kp\,\BBB\leq\al/2$. It follows that $\GG(y)$ satisfies
$$
\e\GG(y(t))\leq\e\GG(y(0))+\kp\,\e\int_0^t \GG(y(s))(-\al\ees(y(s))+C(\gamma,\BBB,\|h\|)\,ds.
$$
It remains to use the inequality
$$
\kp e^v(-\al v+C_1)\leq -\al e^v+C_2\q \text{ for all } v\geq -C,
$$
and the Gronwall lemma, to conclude.
\ep 
\subsection{Exponential supermartingale-type inequality}
The following result provides an estimate for the rate of growth of solutions.
\begin{proposition}\label{supermartingale}
Under the hypotheses of Theorem \ref{2.12}, the following a priori estimate holds for solutions of equation \eqref{1.1}
\be 
\pp \bg\{\sup_{t\geq 0}(\ees_u(t)+\int_0^t(\al \ees_u(s)-\kkk)\,ds)\geq\ees_u(0)+r\bg\}\leq e^{-\beta r} \q\tx{ for any } r>0,
\ee
where $\kkk$ is the constant from inequality \eqref{2.19}, and $\beta=\al/8\cdot(\sup b_j^2)^{-1}$.
\end{proposition}
\bp
Let us first note that
\begin{align*}
\e\sum_{j=1}^\infty b_j^2\int_0^t(\dt u+\al u,e_j)^2\,ds\leq (\sup_{j\geq 1}b_j^2)\int_0^t\e\|\dt u+\al u\|^2\,ds<\infty,\,\quad  t\geq 0.
\end{align*}
It follows that the stochastic integral $M(t)$ defined in \eqref{2.20} is a martingale, and its quadratic variation $\langle M\rangle(t)$ equals 
$$
\langle M\rangle(t)=4\sum_{j=1}^\infty b_j^2\int_0^t(\dt u+\al u,e_j)^2\,ds\leq 4\sup_j b_j^2\int_0^t\|\dt u+\al u\|^2\,ds.
$$
Combining this with inequality inequality \eqref{2.19}, we obtain
$$
\ees_u(t)+\int_0^t(\al \ees_u(s)-\kkk)\,ds\leq\ees_u(0)+\left(M(t)-\f{1}{2}\beta\langle M\rangle(t)\right).
$$
We conclude that
\begin{align*}
&\pp \bg\{\sup_{t\geq 0}(\ees_u(t)+\int_0^t(\al \ees_u(s)-\kkk)\,ds)\geq\ees_u(0)+r)\bg\}\\
&\leq\pp \bg\{\sup_{t\geq 0}(M(t)-\f{1}{2}\beta\langle M\rangle(t))\geq r\bg\}=\pp \bg\{\sup_{t\geq 0}\exp(\beta M(t)-\f{1}{2}\langle \beta M\rangle(t))\geq e^{\beta r}\bg\}\\
&\leq e^{-\beta r},
\end{align*}
where we used the exponential supermartingale inequality.
\ep 

\medskip
\nt
We recall that for a process of the form $y(t)=[u(t),\dt u(t)]$, $\fff^u\equiv\fff_y$ stands for the functional defined by \eqref{2.13}, and $\tau^u\equiv\tau_y$ stands for the stopping time defined by \eqref{6 In-6,1,1}.
\begin{corollary}\label{supermartingale lemma}
Suppose that the hypotheses of Theorem \ref{2.12} are fulfilled. Then for any solution $u(t)$ of equation \eqref{1.1}, we have
\begin{align*}
\pp \bg\{\sup_{t\geq 0}(\fff^u(t)-Lt)\geq\fff^u(0)+r\bg\}&\leq  \exp(4\beta C-\beta r) \q\q\q\q\tx{ for any } r>0,\nt\\
\pp\{l\leq\tau^u<\iin\}&\leq \exp(4\beta C-\beta r-\beta l M)  \q\tx{ for any } l\geq 0,
\end{align*}
where $L=\kkk+4\al C$, $\kkk$ and $\beta$ are the constants from the previous proposition and $C$ is the constant from inequalities \eqref{1.5}-\eqref{1.6}.
\end{corollary}

\nt
This result follows from Proposition \ref{supermartingale} and the fact that, due to inequality \eqref{1.5}, we have
$$
\ees_u(t)\leq|\ees_u(t)|\leq\ees_u(t)+4 C.
$$

\subsection{Existence of stationary measure}
In this subsection we show that the process $y(t)=[u(t),\dt u(t)]$ has a bounded second moment in the more regular space $\h^s=H^{s+1}(D)\times H^s(D)$, with $s=s(\rho)>0$ sufficiently small. By the Bogolyubov-Krylov argument, this immediately implies the existence of stationary distribution for the corresponding Markov process. 
\begin{proposition}\label{2.11}
Under the hypotheses of Theorem \ref{2.12}, there is an increasing function $Q$ such that, for any $s\in(0,1-\rho/2)$, and any solution of equation \eqref{1.1}, we have
$$
\e|y(t)|_{\h^s}^2\leq Q(|y(0)|_{\h})+|y(0)|_{\h^s}^2e^{-\al t}.
$$
\end{proposition}
\bp
Let us split $u$ to the sum $u=v+z$, where $v$ solves
\be \label{2.6}
\p_t^2 v+\gamma \p_t v-\de v=h(x)+\eta(t), \q \xi_v(0)=\xi_u(0).
\ee 
The standard argument shows that for any $s\in[0,1]$, we have
\be 
\e|\xi_v(t)|_{\h^s}^2\leq C(\gamma, \|h\|_1)+|y(0)|_{\h^s}^2 e^{-\al t},
\ee
so that it remains to bound the average of $|\xi_z(t)|_{\h^s}^2$. In view of \eqref{1.1} and \eqref{2.6}, $z(t)$ is the solution of
\be\label{2.7} 
\p_t^2 z+\gamma \p_t z-\de z+f(u)=0,\q \xi_z(0)=0.
\ee
We now follow the argument used in \cite{zelik2004}.
Let us differentiate \eqref{2.7} in time, and set $\theta=\p_t z$. Then $\theta$ solves
\be\label{2.8} 
\p_t^2 \theta+\gamma \p_t \theta-\de \theta+f'(u)\p_t u=0,\q [\theta(0),\dt \theta(0)]=[0,-f(u(0))].
\ee
Let us fix $s\in(0,1-\rho/2)$, multiply this equation by $(-\de)^{s-1}(\dt \theta+\al\theta)$ and integrate over $D$. We obtain
\be\label{2.9}
\f{d}{dt}\tilde\ees_{\theta}(t)+\f{3\al}{2}\tilde\ees_{\theta}(t)\leq 2\int_D |f'(u)\dt u||(-\de)^{s-1}(\dt \theta+\al\theta)|\,dx=:\elll,
\ee
where we set
$$
\tilde\ees_{\theta}(t)=|\xi_\theta|_{\h^{s-1}}^2+\al\gamma|\theta|_{H^{s-1}}^2+2\al(\theta,\dt \theta)_{H^{s-1}}.
$$
By the H\"older and Sobolev inequalities
\begin{align*}
\elll&\leq C_1\int_D (|u|^\rho+1)|\dt u||(-\de)^{s-1}(\dt \theta+\al\theta)|\,dx\\
&\leq C_1(|u|_{L^6}^\rho+1)|\dt u|_{L^2}|(-\de)^{s-1}(\dt \theta+\al\theta)|_{L^{6/(3-\rho)}}\\
&\leq C_2(|u|_{L^6}^2+1)|\dt u|_{L^2}|(-\de)^{s-1}(\dt \theta+\al\theta)|_{H^{1-s}}\\
&\leq C_3(\|\g u\|^2+1)\|\dt u\||\dt \theta+\al\theta|_{H^{1-s}}\leq \f{\al}{2}\tilde\ees_{\theta}(t)+C_4(\|\g u\|^4+1)\|\dt u\|^2,
\end{align*}
where we used the embedding $H^{1-s}\hookrightarrow L^{6/(3-\rho)}$.
Substituting this estimate in \eqref{2.9} and taking the mean value we obtain
$$
\f{d}{dt}\e\tilde\ees_{\theta}(t)\leq -\al\e\tilde\ees_{\theta}(t)+C_4\e(\|\g u\|^4+1)\|\dt u\|^2.
$$
Applying the Gronwall lemma and using Proposition \ref{proposition_exp}, we see that
$$
\e\tilde\ees_{\theta}(t)\leq\e\tilde\ees_{\theta}(0)+C_5,
$$
where the constant $C_5$ depends only on $\al$ and $|y(0)|_\h$.  Moreover, by \eqref{2.8} we have
$$
\tilde\ees_{\theta}(0)=|f(u(0))|_{H^{s-1}}^2\leq |f(u(0))|_{L^2}^2\leq C(1+|y(0)|_\h^6), 
$$
so that 
\be\label{2.10}
\e\tilde\ees_{\theta}(t)\leq Q_1(|y(0)|_\h).
\ee
In view of \eqref{2.7}
\begin{align*}
|z|_{H^{s+1}}=|\de z|_{H^{s-1}}=|\ddot z+\gamma \dt z+f(u)|_{H^{s-1}}&\leq |\ddot z+\gamma \dt z|_{H^{s-1}}+|f(u)|_{L^2}\\
&\leq |\dt \theta+\gamma \theta|_{H^{s-1}}+C(1+|y(0)|_\h^3),
\end{align*}
whence 
$$
|z|_{H^{s+1}}^2\leq 2\tilde\ees_{\theta}(t)+C_6(1+|y(0)|_\h^{6}).
$$
Taking the mean value in this inequality and using \eqref{2.10}, we obtain
$$
\e|\xi_z(t)|_{\h^s}^2\leq Q_2(|y(0)|_\h).
$$
This completes the proof of Proposition \ref{2.11}.
\ep 

\section{Stability of solutions}\label{4.0}
In this section we establish the stability and the recurrence property of solutions of equation \eqref{1.1}.
\subsection{The Foia\c{s}-Prodi estimate}
Here we establish an estimate which will allow us to use the Girsanov theorem.
Let us consider the following two equations:
\begin{align}
\p_t^2 u+\gamma \p_t u-\de u+f(u)&=h(x)+\p_t g(t,x),\label{FP"IN"1}\\
\p_t^2 v+\gamma \p_t v-\de v+f(v)+P_N[f(u)-f(v)]&=h(x)+\p_t g(t,x)\label{FP"IN"2},
\end{align}
where $g(t)$ is a function in $C(\rr_+;H^1_0(D))$, and $P_N$ stands for the orthogonal projection from $L^2(D)$ to  its $N$-dimensional subspace spanned by the functions $e_1,e_2,\ldots,e_N$.

\begin{proposition}\label{4.13}
Suppose that for some non-negative constants $K,l,s$ and $T$ the inequality
\be\label{4.17}
\int_s^t\|\g z\|^2\,d\tau\leq l+K(t-s) \q \text{ for } s\leq t\leq s+T,
\ee
holds for $z=u$ and $z=v$, where $u$ and $v$ are solutions of \eqref{FP"IN"1} and \eqref{FP"IN"2}, respectively. Then, for any $\es>0$ there is an integer $N_*\geq 1$ depending only on $\es$ and $K$ such that for all $N\geq N_*$ we have
\be \label{4.16}
|\xi_v(t)-\xi_u(t)|^2_\h\leq e^{-\al(t-s)+\es l}|\xi_v(s)-\xi_u(s)|^2_\h \q \text{ for } s\leq t\leq s+T.
\ee
\end{proposition}

\bp 
 Let us set $w=v-u$. Then $w(t)$ solves
\be\label{6.50}
\p_t^2 w+\gamma \p_t w-\de w+(I-P_N)[f(v)-f(u)]=0,
\ee
and we need to show that the flow $y(t)=\xi_w(t)$ satisfies
\be\label{4.15}
|y(t)|^2_\h\leq e^{-\al (t-s)+\es l}|y(s)|^2_\h \q \text{ for } s\leq t\leq s+T.
\ee
The function $y(t)$ satisfies
\begin{align}
\p_t|y|_{\h}^2&=2[(\g w,\g \dt w)+(-\gamma \dt w+\de w-(I-P_N)[f(v)-f(u)]+\al \dt w,\dt w+\al w)]\notag\\
&\leq 2[(\g w,\g \dt w)
+(-\gamma \dt w+\de w+\al \dt w,\dt w+\al w)]\notag\\
&\q+2\|(I-P_N)[f(v)-f(u)]\|(\|\dt w\|+\al\|w\|)\notag\\
&\leq -\f{3\al}{2} |y|_\h^2+4\|(I-P_N)[f(v)-f(u)]\||y|_\h\label{4.8}.
\end{align}
We first note that
\be\label{4.7}
\|(I-P_N)[f(v)-f(u)]\|\leq|I-P_N|_{L(H^{1,p}\to L^2)}|f(v)-f(u)|_{H^{1,p}},
\ee
and that
\be\label{FP"IN"3}
|f(v)-f(u)|_{H^{1,p}}\leq C\sum_{j=1}^3|\p_j[f(v)-f(u)]|_{L^p},
\ee
where $p\in(6/5,2)$ will be chosen later.
Further,
\begin{align}
|\p_j[f(v)-f(u)]|_{L^p}&=|f'(v)\p_j v-f'(u)\p_j u|_{L^p}\notag\\
&\leq |(f'(v)-f'(u))\p_j v|_{L^p}+|f'(u)(\p_j v-\p_j u)|_{L^p}=J_1+J_2\label{FP"IN"4}.
\end{align}
For $J_1$ we have 
\begin{align}
J_1&=\left(\int_D |(f'(v)-f'(u))\p_j v|^p\,dx\right)^{\f{1}{p}}\notag\\
&\leq C_4\left(\int_D |w|^p|\p_j v|^p(|v|^{p(\rho-1)}+|u|^{p(\rho-1)}+1)\,dx\right)^{\f{1}{p}}\notag\\
&\leq C_5|w|_{L^6}|\g v|_{L^2}(|v|_{L^6}^{\rho-1}+|u|_{L^6}^{\rho-1}+1)\leq C_6\|w\|_1(\|v\|_1^2+\|u\|_1^2+1)\label{FP"IN"5},
\end{align}
where we used the H\"older and Sobolev inequalities and chose $p=6(3+\rho)^{-1}$.
And finally, for $J_2$ we have
\begin{align}
J_2=\left(\int|f'(u)|^p|\p_j v-\p_j u|^p\,dx\right)^{\f{1}{p}}&\leq C_7\left(\int|\p_j w|^p (|u|^{p\rho}+1)\,dx\right)^{\f{1}{p}}\notag\\
&\leq C_8\|w\|_1(\|v\|_1^2+\|u\|_1^2+1)\label{FP"IN"6},
\end{align}
where we once again used the H\"older inequality.
Combining inequalities \eqref{4.7}-\eqref{FP"IN"6} together, we obtain
\be \label{4.10}
\|(I-P_N)[f(v)-f(u)]\|\leq C'_1|I-P_N|_{L(H^{1,p}\to L^2)}(\|v\|_1^2+\|u\|_1^2+1)|y|_\h.
\ee
Substituting this inequality in \eqref{4.8}, we see that
\be\label{4.14} 
\p_t|y|_{\h}^2\leq (-3\al/2+C|I-P_N|_{L(H^{1,p}\to L^2)}(\|v\|_1^2+\|u\|_1^2+1))|y|_\h^2.
\ee
By the Sobolev embedding theorem, the space $H^{1,p}(D)$ is compactly embedded in $L^2(D)$ for $p>6/5$. This implies that the sequence $|I-P_N|_{L(H^{1,p}\to L^2)}$ goes to zero as $N$ goes to infinity. Combining this fact with the Gronwall lemma applied to \eqref{4.14} and using \eqref{4.17}, we arrive at \eqref{4.15}.
\ep
\subsection{Controlling the growth of intermediate process}
The goal of this subsection is to show that inequality \eqref{4.17} (and therefore \eqref{4.16}) holds with high probability, for $g(t)=\zeta(t)$.

\nt
For any $\h$-valued continuous process $y(t)$, let $\tau_y$ be the stopping time defined in \eqref{6 In-6,1,1},
where $L$ is the constant constructed in Corollary \ref{supermartingale lemma}, and $M,r$ are some positive constants. We recall that for the process of the form $y=[z,\dt z]$ we shall write $\tau^z$ instead of $\tau_{[z,\dt z]}$.
\begin{proposition}\label{4.23}
Let $u$ and $v$ be solutions of \eqref{FP"IN"1} and \eqref{FP"IN"2} where $g(t)=\zeta(t)$, that are issued from initial points $y,y'\in B_1$, respectively. Then
\be \label{4.18}
\pp\{\tau^v<\iin\}\leq 3 \exp(4\beta C-\beta r)+C_{M,r}|y-y'|_\h,
\ee
where $\beta$ is the constant from Proposition \ref{supermartingale}.
\end{proposition}
\bp 
To prove this result, we follow the arguments presented in Section 3.3 of \cite{KS-book} and Section 4 of \cite{kuksin2013stochastic}. First, note that since inequality \eqref{4.18} concerns only the law of $v$ and not the solution itself, we are free to choose the underlying probability space $(\omm,\fff,\pp)$. We assume that it coincides with the canonical space of the Wiener process $\{\hat\zeta (t)\}_{t\geq 0}$. More precisely, $\omm$ is the space of continuous functions $\om:\rr_+\to\h$ endowed with the metric of uniform convergence on bounded intervals, $\pp$ is the law of $\hat\zeta$ and $\fff$ is the completion of the Borel $\sigma$-algebra with respect to $\pp$.

\nt
Let us define vectors $\hat e_j=[0,e_j]$ and their vector span
$$
 \h_N=\text{span}\{\hat e_1,\hat e_2,\ldots,\hat e_N\},
 $$
 which is an $N$-dimensional subspace of $\h$. The space $\Omega=C(\rr_+,\h)$ can be represented in the form
 $$
 \Omega=\Omega_N\dt+\Omega_N^{\perp},
 $$
 where $\Omega_N=C(\rr_+, \h_N)$ and $\Omega_N^\perp=C(\rr_+, {\h^\perp_N})$.
We shall write $\om=(\om^{(1)},\om^{(2)})$ for $\om=\om^{(1)}\dt+\om^{(2)}$.

Let $u'$ be a solution of equation \eqref{FP"IN"1} that has the same initial data as $v$. Introduce the stopping time
\be \label{4.24}
\tilde\tau=\tau^u\wedge\tau^{u'}\wedge\tau^v,
\ee
and a transformation $\Phi:\Omega\to\Omega$ given by
\be\label{4.22}
 \Phi(\om)(t)=\om(t)-\int_0^t a(s)\,ds,\q a(t)=\ch_{t\leq\tilde\tau}\PP_N(0,[f(u)-f(v)]),
\ee
 where $\PP_N$ is the orthogonal projection from $\h$ to $\h_N$.
\begin{lemma}\label{4.21}
For any initial points
 $y$ and $y'$ in $B_1$, we have
\be \label{4.19}
|\pp-\Phi_*\pp|_{var}\leq C_{M,r} |y-y'|_\h,
\ee
where $\Phi_*\pp$ stands for the image of $\pp$ under $\Phi$.
\end{lemma}

\nt
{\it Proof of lemma \ref{4.21}}.

\nt
{\it Step 1.} 
Let us note that by the definition of $\tilde\tau$ we have
\be\label{6.54}
\fff^u(t)\leq \fff^u(0)+(L+M)t+r,\q
\fff^v(t)\leq\fff^{u'}(0)+(L+M)t+r,
\ee
for all $t\leq\tilde\tau$.
We claim that there is an integer $N=N(\al,L,M)$ such that for all $t\leq\tilde\tau$ we have
\be\label{7.8}
|\xi_{v}(t)-\xi_{u}(t)|_\h^2\leq e^{-\al t+\theta}|\xi_u(0)-\xi_{u'}(0)|_\h^2,
\ee
where $\theta=|\ees_u(0)|\vee|\ees_{u'}(0)|+r$. Indeed, in view of inequality \eqref{1.5}, for any $y=[y_1,y_2]$ in $\h$, we have
\begin{align*}
|y(t)|_{\h}^2&\leq |y(t)|_{\h}^2+2\int_D F(y_1)\,dx|-2\int_D F(y_1)\,dx
\leq |\ees(y)|+2\nu \|y_1\|^2+2C\\
&\leq |\ees(y)|+\f{\lm_1}{2}\|y_1\|^2+2C\leq |\ees(y)|+\f{1}{2}|y|_\h^2+2C,
\end{align*}
so that
\be\label{7.9}
|y|_{\h}^2\leq 2 |\ees(y)|+4C.
\ee
Combining this inequality with \eqref{6.54}, we see that for all $t\leq\tilde\tau$
$$
\al\int_0^t\|\g z(s)\|^2\,ds \leq 2(|\ees_u(0)|\vee |\ees_u'(0)|+r)+2(L+M+2C)t,
$$
for $z=u$ and $z=v$. Using this inequality and applying Proposition \ref{4.13} with $\es=\al/2$ we arrive at \eqref{7.8}.

\nt
{\it Step 2.} 
Let us note that the transformation  $\Phi$ can be represented in the form
$$
\Phi(\om)=(\Psi(\om),\om^{(2)}),
$$
where $\Psi:\omm\to\omm_N$ is given by
$$
\Psi(\om)(t)=\om^{(1)}(t)+\int_0^t a(s;\om)\,ds.
$$
It is straightforward to see that 
 \be \label{4.20}
 |\pp-\Phi_*\pp|_{var}\leq\int_{\omm^\perp_N}|\Psi_*(\pp_N,\om^{(2)})-\pp_N|_{var}\pp^\perp_N(d\om^{(2)}),
 \ee
 where $\pp_N$ and $\pp^\perp_N$ are the images of $\pp$ under the natural projections $\PP_N:\omm\to\omm_N$ and $\QQ_N:\omm\to\omm^\perp_N$, respectively.
Define the processes 
$$
 z(t)=\om^{(1)}(t),\q
 \tilde z(t)=\om^{(1)}(t)+\int_0^t a(s;\om)\,ds.
$$
 It follows that $\pp_N=\dd z$ and $\Psi_*(\pp,\om^{(2)})=\dd\tilde z$.
 By Theorem A.10.1 in \cite{KS-book}, we have
 \be\label{Dz-D tilde z}
 |\dd z-\dd \tilde z|_{var}\leq\f{1}{2}\left(\left(\e\exp\left[6\max_{1\leq j\leq N} b_j^{-1}\int_0^\infty |a(t)|^2\,dt\right]\right)^{\f{1}{2}}-1\right)^{\f{1}{2}},
 \ee
 provided the Novikov condition
 $$
 \e\exp\left(C\int_0^\infty |a(t)|^2\,dt\right)<\infty,\q \text{ for any } C>0,
 $$
 holds. In view of inequalities \eqref{6.54} and \eqref{7.8} we have
\begin{align*}
&\e\exp(C\int_0^\infty |a(t)|^2\,dt)=\e\exp(C\int_0^{\tilde\tau} |a(t)|^2\,dt)\\
&\leq\e\exp(C\int_0^{\tilde\tau} \|f(v)-f(u)\|^2\,dt)\leq\e\exp(C_1\int_0^{\tilde\tau} \|v-u\|_1^2(1+\|u\|_1^4+\|v\|_1^4)\,dt)\\
&\leq\e\exp(C_2 |\xi_u(0)-\xi_{u'}(0)|_\h^2\int_0^\infty e^{-\al t+\theta}K(t)\,dt),
 \end{align*}
 where
 $$
 K(t)=(1+|\ees_u(0)|\vee|\ees_{u'}(0)|+(L+M)t+r)^2.
 $$
 So not only the Novikov condition holds, but also there is a positive constant $C_{M,r}=C(\al,L,M,r)$ such that the term on the right-hand side of inequality \eqref{Dz-D tilde z} does not exceed $C_{M,r} |y-y'|_\h$. Combining this with inequality \eqref{4.20}, we arrive at \eqref{4.19}.
 
 \nt
Now we are ready to establish \eqref{4.18}. Introduce auxiliary $\h$-continuous processes $y^u, y^{u'}$ and $y^v$ defined as follows: for $t\leq \tilde\tau$ they coincide with processes $\xi_u, \xi_{u'}$ and $\xi_v$, respectively, while for $t\geq\tilde\tau$ they solve 
$$ 
\p_t y=-\lm y,
$$
where $\lm>0$ is a large parameter.
By construction, with probability 1, we have
\be\label{4.29}
y^v(t,\om)=y^{u'}(t,\Phi(\om)) \q\text{ for all } t\geq 0.
\ee
Let us note that
\begin{align}\label{4.25}
\pp(\tau^v<\iin)&=\pp(\tau^v<\iin, \tau^u\wedge\tau^{u'}<\iin)+\pp(\tau^v<\iin, \tau^u\wedge \tau^{u'}=\iin)\notag\\
&\leq \pp(\tau^u<\iin)+\pp(\tau^{u'}<\iin)+\pp(\tau^v<\iin, \tau^u\wedge \tau^{u'}=\iin).
\end{align}
Moreover, in view of \eqref{4.29}
\begin{align}
&\pp(\tau^v<\iin, \tau^u\wedge\tau^{u'}=\iin)\leq\pp(\tau_{y^v}<\iin)=\Phi_*\pp(\tau_{y^{u'}}<\iin)\notag\\
&\q\leq \pp(\tau_{y^{u'}}<\iin)+|\pp-\Phi_*\pp|_{var}\leq \pp(\tau^{u'}<\iin)+|\pp-\Phi_*\pp|_{var}\label{6.43},
\end{align}
where we used the fact that for $t\geq \tilde\tau$ the norms of auxiliary processes decay exponentially. Combining these two inequalities we obtain
\begin{align}\label{4.28}
\pp(\tau^v<\iin)\leq \pp(\tau^u<\iin)+2\pp(\tau^{u'}<\iin)+|\pp-\Phi_*(\pp)|_{var}.
\end{align}
It remains to use Corollary \ref{supermartingale lemma} and Lemma \ref{4.21}  to conclude.
\ep

\subsection{Hitting a non-degenerate ball}
Here we show that the trajectory of the process $y(t)=[u(t),\dt u(t)]$ issued from arbitrarily large ball hits any non-degenerate ball centered at the origin, with positive probability, at a finite non-random time. We denote by $B_d$ the ball of radius $d$ in $\h$, centered at the origin.
\begin{proposition}\label{6.12}
For any $R>0$ and $d>0$ there is $T_*=T_*(R,d)>0$ such that for all $T\geq T_*$, we have
\be\label{6.13}
\inf_{y\in B_R}P_T(y, B_d)>0,
\ee
where $P_t(y,\Gamma)=\pp_y(S_t(y,\cdot)\in\Gamma)$ is the transition function of the Markov process corresponding to \eqref{1.1}.
\end{proposition}
\bp 
Let us first split $u$ to the sum $\tilde u+\bar u$, where $\bar u$ is the solution of
\be\label{6.15} 
\p_t^2 \bar u+\gamma\p_t \bar u-\de \bar u=0, \q\xi_{\bar u}(0)=\xi_u(0).
\ee
Then the corresponding flow $\bar y(t)$ satisfies the exponential decay estimate
\be\label{6.14}
|\bar y(t)|_\h^2\leq e^{-\al t}|y(0)|_\h^2.
\ee
Let us fix $T_*=T_*(R,d)$ such that for all $T\geq T_*$ and any initial point $y$ in $B_R$, we have
\be\label{6.23} 
|\bar y(T)|_\h\leq d/2.
\ee
We claim that \eqref{6.13} holds with this time $T_*$. Indeed, if this is not true, then there is $T\geq T_*$ such that 
\be\label{6.58}
\inf_{y\in B_R}P_T(y, B_d)=0.
\ee
In view of \eqref{6.15}, $\tilde u$ solves
\be\label{6.16}
\p_t^2 \tilde u+\gamma \p_t \tilde u-\de \tilde u+f(\tilde u+\bar u)=h(x)+\eta(t,x),\q \xi_{\tilde u}(0)=0.
\ee 
Now note that this equation is equivalent to 
\be
\tilde y(t)=\int_0^t g(s)\,ds+\hat\zeta(t), 
\ee
where
$$
\tilde y=[\tilde u,\p_t \tilde u],\q
g=[\p_t \tilde u,-\gamma\p_t \tilde u+\de \tilde u-f(\tilde u+\bar u)+h(x)],\q
\hat \zeta=[0,\zeta],
$$
and therefore $\tilde y(T)$ continuously depends on $\hat \zeta$ (in the sense that the small perturbation of $\hat\zeta$ in ${C(0,T;\h)}$ will result in a small perturbation of $\tilde y(T)$ in $\h$). Let us consider equation \eqref{6.16} with the right-hand side
\be \label{6.24}
\tilde\zeta^y=\int_0^t f(\bar u)\,ds-th,
\ee
which is a non-random force (the notation $\tilde\zeta^y$ is justified by the fact that it is uniquely determined by the initial point $y$). 
Then the function $\tilde u\equiv 0$ solves that equation. It follows that there exists $\es=\es(d)>0$ such that
$$
|\tilde y(T)|_\h\leq d/2,
$$
provided
\be
|\zeta-\tilde\zeta^y|_{C(0,T;L^2)}\leq\es. 
\ee
Combining this with inequality \eqref{6.23} we obtain 
\be 
|y(T)|_\h\leq d.
\ee
Therefore
\be \label{6.18}
P_T(y,B_d)\geq \pp(|\zeta-\tilde\zeta^y|_{C(0,T;L^2)}\leq\es).
\ee
We need the following lemma. It is established in the appendix.
\begin{lemma}\label{6.20}
For any $\rho<2$ there exists $s=s(\rho)>0$ such that if 
\be \label{5.17}
|f'(u)|\leq C(|u|^{\rho}+1),
\ee
then
\be\label{6.39}
|f(u)-f(v)|_{L^2}\leq C_1(|u|_{H^{1-s}}^\rho+|v|_{H^{1-s}}^\rho+1)|u-v|_{H^{1-s}},
\ee
where $C_1>0$ depends only on $C>0$.
\end{lemma}
Let us suppose that we have \eqref{6.58}, and let $y_j(0)=[u_j(0),\dt u_j(0)]$ be a minimizing sequence. This sequence is bounded in $H^1\times L^2$, so it has a converging subsequence in $H^{1-s}\times H^{-s}$ ($s$ is the constant from the previous lemma). Moreover, a standard argument coming from theory of $m$-dissipative operators shows that the resolving operator of \eqref{6.15} generates a continuous semigroup in $H^{1-s}\times H^{-s}$ (e.g., see \cite{cazenave1998introduction}). It follows that for all $t\geq 0$ the corresponding sequence of solutions $\bar y_j(t)$ issued from $y_j(0)$ converges in that space. In particular $\bar u_j(t)$ converges in $H^{1-s}$. Denoting by $\hat u(t)$ its limit and using Lemma \ref{6.20} together with inequality 
$$ 
|\int_0^t \psi(s)\,ds|_{C(0,T;L^2)}\leq T^{1/2}|\psi|_{L^2(0,T;L^2)},
$$
we see that $\tilde\zeta_j=\tilde\zeta^{y_j(0)}\to\tilde\zeta$ in $C(0,T;L^2)$, where
$$
\tilde\zeta=\int_0^t f(\hat u)\,ds-th.
$$
Inequality \eqref{6.18} implies that
\be \label{6.22}
\pp(|\zeta-\tilde\zeta_j|_{C(0,T;L^2)}\leq\es)\to 0 \q \text{ as } j\to\iin.
\ee
Let us fix $j_0\geq 1$ so large that for all $j\geq j_0$
$$
|\tilde\zeta_j-\tilde\zeta|_{C(0,T;L^2)}\leq \es/2.
$$
Then by the triangle inequality, for all $j\geq j_0$
\begin{align*}
\pp(|\zeta-\tilde\zeta|_{C(0,T;L^2)}\leq \es/2)&=\pp(|\zeta-\tilde\zeta_j+\tilde\zeta_j-\tilde\zeta|_{C(0,T;L^2)}\leq\es/2)\\
&\leq \pp(|\zeta-\tilde\zeta_j|_{C(0,T;L^2)}\leq |\tilde\zeta_j-\tilde\zeta|_{C(0,T;L^2)}+\es/2)\\
&\leq \pp(|\zeta-\tilde\zeta_j|_{C(0,T;L^2)}\leq \es).
\end{align*}
Letting $j$ go to $\iin$ and using inequality \eqref{6.22}, we obtain
$$
\pp(|\zeta-\tilde\zeta|_{C(0,T;L^2)}\leq \es/2)=0,
$$
which is impossible, since the support of $\zeta$ restricted to $[0,T]$ coincides with ${C(0,T;L^2)}$. The proof of Proposition \ref{6.12} is complete.
\ep

\section{Proof of Theorem \ref{2.14}}\label{6.31}
In this section we establish Theorem \ref{2.14}. As it was already mentioned, this will imply Theorem \ref{2.15}. We then show that the non-degeneracy condition imposed
on the force can be relaxed to allow forces that are non-degenerate only in the
low Fourier modes (see Theorem \ref{6.10}).

\medskip
\subsection{Recurrence: verification of \eqref{2.17}-\eqref{2.16}}\label{6.32}
In view of Proposition \ref{proposition_exp}, it is sufficient to establish inequality \eqref{2.16}. To this end, we shall use the existence of a Lyapunov function, combined with an auxiliary result established in \cite{shirikyan-bf2008}.

Let $S_t(y,\om)$ be a Markov process in a separable Banach space $\h$ and let $\RR_t(\yy,\om)$ be its extension on an interval $[0, T]$. Consider a continuous functional $\gi(y)\geq 1$ on $\h$ such that
$$
\lim_{|y|_\h\to\infty}\gi(y)=\infty.
$$ 
Suppose that there are positive constants $d, R, t_*, C_*$ and $a<1$, such that
\begin{align}
\e_y\gi(S_{t_*})&\leq a\, \gi(y)\,\,\q\tx{ for } |y|_\h\geq R,\label{6.55}\\
\e_y\gi(S_{t})&\leq C_*\q\q\,\,\,\,\,\tx{ for } |y|_\h\leq R,\, t\geq 0,\label{6.56}
\end{align}
\be
\inf_{y,y'\in B_R}\pp_\yy\{|\R_T(y,y',\cdot)|_{\h}\vee |\R'_T(y,y',\cdot)|_{\h}\leq d\}>0.\label{6.57}
\ee

\medskip\nt
We shall denote by $\tau_d$ the first hitting time of the set $B_{\hh}(d)$. The following proposition is a weaker version of the result proved in \cite{shirikyan-bf2008} (see Proposition 3.3).
\begin{proposition}\label{6. Proposition-Lyapunov-dissipation}
Under the above hypotheses there are positive constants $C$ and $\kp$ such that the inequality
\be 
\e_\yy\exp(\kp\tau_d)\leq C(\gi(y)+\gi(y')), \q\tx{ for any } \yy=(y,y')\in\hh,
\ee
holds for the extension $\SSS_t$ constructed by iteration of $\RR_t$ on the half-line $t\geq 0$.
\end{proposition}

\nt
It follows from estimate \eqref{2.2} that inequalities \eqref{6.55} and \eqref{6.56} are satisfied for the functional
$$
\gi(y)=1+|\ees(y)|.
$$
We now show that for any $d>0$ we can find an integer $k\geq 1$ and $T_*\geq 1$ sufficiently large, such that we have \eqref{6.57} for any $T\in\{kT_*,(k+1)T_*,\ldots\}$. In what follows, we shall drop the subscript and write $|y|$ instead of $|y|_\h$. So let us fix any $d>0$, and consider the events
\begin{align*}
G_d&=\{|\R_T(y,y')|\leq d\},\q G'_d=\{|\R'_T(y,y')|\leq d\},\\
E_r&=\{\fff_\R(t)\leq \fff_\R(0)+Lt+r\}\cap \{\fff_\R'(t)\leq \fff_\R'(0)+Lt+r\},
\end{align*}
where $\fff_y(t)$ is defined in \eqref{2.13}, and $L$ is the constant from Corollary \ref{supermartingale lemma}.

\nt
{\it Step 1.} 
First, let us note that by Proposition \ref{6.12}, there is $T_*=T_*(R,d)\geq 1$, such that
\be\label{7.5}
\pp_y\{|S_{T_*}(y,\cdot)|\leq d/2\}\geq c_d \q\text{ for any } y\in B_R,
\ee
where $c_d$ is a positive constant depending on $d,R$ and $T_*$.
We claim that this implies 
\be\label{7.3}
\pp_y\{|S_{kT_*}(y,\cdot)|\leq d/2\}\geq c_d \q\text{ for all } k\geq 1 \text{ and } y\in B_R.
\ee
Indeed, let us fix any integer $k\geq 1$ and introduce the stopping times
$$
\bar\tau(y)=\min\{nT_*, n\geq 1 : |S_{nT_*}(y,\cdot)|>d/2\},\q \bar\sigma=\bar\tau\wedge kT_*.
$$
Let us note that if $\bar\tau$ is finite, then we have
\be\label{7.6}
|S_{\bar\tau-T_*}(y,\cdot)|\leq R \q\text{ and }\q |S_{\bar\tau}(y,\cdot)|>d/2,
\ee
where inequalities hold for any $y$ in $B_R$.
Moreover
\be\label{7.4}
\pp_y\{|S_{kT_*}(y,\cdot)|> d/2\}\leq \pp_y\{\bar\tau=\bar\sigma\},
\ee
where we used that for $\bar\tau>kT_*$, we have $|S_{kT_*}(y,\cdot)|{_\h}\leq d/2$. In view of \eqref{7.6}
\be\label{7.7}
\pp_y\{\bar\tau=\bar\sigma\}\leq \pp_y\{|S_{\bar\sigma-T_*}(y,\cdot)|\leq R,\,|S_{\bar\sigma}(y,\cdot)|>d/2\}:=p.
\ee
Since $\bar\sigma$ is a.s. finite, we can use the strong Markov property, and obtain
\begin{align*}
p&=\e_y[\e_y(\ch_{|S_{\bar\sigma-T_*}(y,\cdot)|\leq R}\cdot\ch_{|S_{\bar\sigma}(y,\cdot)|>d/2}|\fff_{\bar\sigma-T_*})]\\
&=\e_y[\ch_{|v|\leq R}\cdot\e_v(\ch_{|S_{T_*}(v,\cdot)|>d/2})]\\
&=\e_y[\ch_{|v|\leq R}\cdot\pp_v(|S_{T_*}(v,\cdot)|>d/2)]\leq \sup_{\bar v\in B_R}\pp_{\bar v}(|S_{T_*}(\bar v,\cdot)|>d/2),
\end{align*}
where $v=S_{\bar\sigma-T_*}(y,\cdot)$, and $\fff_t$ is the filtration generated by $S_t$. In view of \eqref{7.5}, the last term in this inequality does not exceed $1-c_d$. Combining this with inequalities \eqref{7.4} and \eqref{7.7}, we arrive at \eqref{7.3}.

\nt
{\it Step 2.} 
It follows from the previous step that for any $T\in\{T_*, 2T_*, \ldots\}$
\be\label{6.59}
\pp_\yy(G_{d/2})\wedge \pp_\yy(G'_{d/2})\geq c_d,\q \text{ for any } y,y'\in B_R,
\ee
where we used that $\RR_t$ is an extension of $S_t$.
Further, by Corollary \ref{supermartingale lemma} we have
\be\label{6.62}
\pp_\yy(E_r)\geq 1-2\exp(4\beta C-\beta r):=1-o(r).
\ee
Let us fix $r=r(d,R,T_*)>0$ so large that
\be\label{6.60}
o(r)\leq c_d^2/8.
\ee
By the symmetry, we can assume that
\be\label{6.63}
\pp_\yy(G'_{d/2}\nnn^c)\leq\pp_\yy(G_{d/2}\nnn^c),
\ee
where we set $\nnn=\{\vvv(y,y')\neq \vvv'(y,y')\}$. We claim that
\be\label{inclusion G-Delta}
G_{d/2}E_r\,\nnn^c\subset G_d G'_d,
\ee
for any $T\in\{kT_*,(k+1)T_*,\ldots\}$ with $k\geq 1$ sufficiently large.
To prove this, let us fix any $\om$ in $G_{d/2}E_r\,\nnn^c$, and note that it is sufficient to establish
\be\label{6.49}
|\R_T(y,y',\om)-\R'_T(y,y',\om)|\leq d/2,\q\text{ for any } y,y' \text{ in } B_R.
\ee
Since $\om\in\nnn^c$, we have that $\vvv=\vvv'$, and therefore, in view of \eqref{7.1}-\eqref{6. Coupling relation 1}, $\R_t(y,y')$ and $\R'_t(y,y')$ are, respectively, the flows of equations
\be 
\p_t^2 \tilde u+\gamma \p_t \tilde u-\de\tilde u+f(\tilde u)-P_N f(\tilde u)=h(x)+\psi(t),\q \xi_{\tilde u}(0)=y,
\ee
and 
\be 
\p_t^2 \tilde v+\gamma \p_t \tilde v-\de\tilde v+f(\tilde v)-P_N f(\tilde v)=h(x)+\psi(t),\q \xi_{\tilde u}(0)=y'.
\ee
It follows that their difference $w=\tilde v-\tilde u$ solves 
$$
\p_t^2 w+\gamma \p_t w-\de w+(I-P_N)[f(\tilde v)-f(\tilde u)]=0,\q [w(0),\dt w(0)]=y'-y.
$$
Using the Foia\c{s}-Prodi estimate established in Proposition \ref{4.13} (see \eqref{6.50}-\eqref{4.15}) together with the fact that $\om\in E_r$, we can find an integer $N\geq 1$ depending only on $L$ such that 
$$
 |\R_T(y,y',\om)-\R'_T(y,y',\om)|^2\leq C(r,R)e^{-\al T}|y-y'|^2\leq 4R^2 C(r,R)e^{-\al T}.
 $$
 Since $r$ is fixed, we can find $k\geq 1$ sufficiently large, such that the right-hand side of this inequality is less than $d^2/4$ for any $T\in\{kT_*,(k+1)T_*,\ldots\}$, so that we have \eqref{inclusion G-Delta}.
 
 \nt
 {\it Step 3.} 
We now follow the argument presented in \cite{shirikyan-bf2008}. In view of \eqref{inclusion G-Delta}
\begin{align*}
\pp_\yy(G_d G'_d)&=\pp_\yy(G_d G'_d \nnn^c)+\pp_\yy(G_d G'_d \nnn)\\
&\geq \pp_\yy(G_d G'_d E_r\nnn^c)+\pp_\yy(G_d|\nnn) \pp_\yy(G'_d|\nnn)\pp_\yy(\nnn)\\
&\geq \pp_\yy(G_{d/2} E_r \nnn^c)+\pp_\yy(G_d \nnn)\pp_\yy(G'_d\nnn),
\end{align*}
where we used the independence of $\vvv$ and $\vvv'$ conditioned on the event $\nnn$. Combining this inequality with \eqref{6.62}, we obtain
$$
\pp_\yy(G_d G'_d)\geq\pp_\yy(G_{d/2}\nnn^c)+\pp_\yy(G_{d} \nnn)\pp_\yy(G'_{d} \nnn)-o(r).
$$
We claim that the right-hand side of this inequality is no less than $c_d^2/8$.
Indeed, if $\pp_\yy(G_{d/2}\nnn^c)\geq c_d^2/4$, then the required result follows from inequality \eqref{6.60}. If not, then by inequalities \eqref{6.59} and \eqref{6.63}, we have
$$
c_d^2\leq \pp_\yy(G_{d/2})\pp_\yy(G'_{d/2})\leq \pp_\yy(G_{d/2}\nnn)\pp_\yy(G'_{d/2}\nnn)+3c_d^2/4,
$$
so that
$$
\pp_\yy(G_{d}\nnn)\pp_\yy(G'_{d}\nnn)\geq\pp_\yy(G_{d/2}\nnn)\pp_\yy(G'_{d/2}\nnn)\geq c_d^2/4.
$$
We have thus shown that for any $y,y'$ in $B_R$
$$
\pp_\yy\{|\R_T(y,y',\cdot)|\vee |\R'_T(y,y',\cdot)|\leq d\}\equiv \pp_\yy(G_d G'_d)\geq c_d^2/8,
$$
and therefore we have \eqref{6.57}. The hypotheses of Proposition \ref{6. Proposition-Lyapunov-dissipation} are thus satisfied, so that inequality \eqref{2.16} holds.\subsection{Exponential squeezing: verification of \eqref{6.33}-\eqref{6.5}}
Let $u,u',v,\tilde u, \tilde u', \tilde v$ and $\vo,\tau,\sigma$ be the processes and stopping times constructed in Subsection \ref{Main result and scheme of its proof}. Consider the following events:
$$
\qqq'_k=\{kT\leq\tau\leq (k+1)T,\tau\leq\vo\},\q\qqq''_k=\{kT\leq\vo\leq (k+1)T,\vo<\tau\}.
$$
\begin{lemma}\label{6.27}
There exist positive constants $d,r, L$ and $M$ such that for any initial point $\yy\in B_{\hh}(d)$ and any $T\geq1$ sufficiently large 
$$
\pp_{\yy}(\qqq'_k)\vee\pp_\yy(\qqq_k'')\leq e^{-2(k+1)}\q\text{ for all }k\geq 0.
$$
\end{lemma}
\bp

\nt\\
{\it~Step1.} (Probability of $\qqq_k'$).
Let $L$ be the constant from Corollary \ref{supermartingale lemma}. Then using second inequality of this corollary, we obtain
\be\label{6.48}
\pp_\yy(\qqq'_k)\leq\pp_\yy(kT\leq\tau<\iin)\leq 2\exp(4\beta C-\beta r-\beta k T M)\leq e^{-2(k+2)},
\ee
for $M\geq 2\beta^{-1}, r\geq 5\beta^{-1}+4C$. From now on, the constants $L,M$ and $r$ will be fixed.

\smallskip\nt
{\it~Step2.} (Probability of $\qqq_k''$). Let us first note that by the Markov property we have
\begin{align} 
\pp_\yy(\qqq''_k)&=\pp_\yy(\qqq''_k,\sigma\geq kT)=\e_\yy(\ch_{\qqq''_k}\cdot\ch_{\sigma\geq kT})=\e_\yy[\e_\yy(\ch_{\qqq''_k}\cdot\ch_{\sigma\geq kT}|\FFF_{kT})]\notag\\
&=\e_\yy[\ch_{\sigma\geq kT}\cdot\e_\yy(\ch_{\qqq''_k}|\FFF_{kT})]\leq \e_\yy[\ch_{\sigma\geq kT}\cdot\e_{\bar\yy}\ch_{0\leq \vo\leq T}]\notag\\
&=\e_\yy[\ch_{\sigma\geq kT}\cdot\pp_{\bar \yy}(0\leq \vo\leq T)],\label{6.28}
\end{align}
where $\bar\yy(\cdot)=\yy(kT,\cdot)$, and $\FFF_t$ stands for the filtration corresponding to the process $\SSS_t$.
Moreover, it follows from the definition of maximal coupling, that for any $\yy$ in $\hh$, we have
$$
\pp_{\yy}(0\leq \vo\leq T)= |\pp_{\yy}\{\xi_v\}_T-\pp_{\yy}\{\xi_{u'}\}_T|_{var}.
$$
Combining this with inequality \eqref{6.28}, we obtain 
\be\label{6.45}
\pp_\yy(\qqq''_k)\leq \e_\yy(\ch_{\sigma\geq kT}\cdot|\pp_{\bar\yy}\{\xi_v\}_T-\pp_{\bar\yy}\{\xi_{u'}\}_T|_{var})
\ee
Further, let us note that
\begin{align}
&|\pp_{\bar\yy}\{\xi_v\}_T-\pp_{\bar\yy}\{\xi_{u'}\}_T|_{var}=\sup_{\Gamma}|\pp_{\bar\yy}(\{\xi_v\}_T\in\Gamma)-\pp_{\bar\yy}(\{\xi_{u'}\}_T\in\Gamma)\notag|\\
&\q\leq \pp_{\bar\yy}(\tilde\tau<\iin)+ \sup_{\Gamma}|\pp_{\bar\yy}(\{\xi_v\}_T\in\Gamma, \tilde\tau=\iin)
-\pp_{\bar\yy}(\{\xi_{u'}\}_T\in\Gamma, \tilde\tau=\iin)|\notag\\
&\q:=\elll_1+\elll_2\label{6.41},
\end{align}
where $\tilde\tau=\tau^u\wedge\tau^{u'}\wedge\tau^v$, and the supremum is taken over all $\Gamma\in\bbb(C(0,T;\h))$. In view of \eqref{4.29} we have
\be\label{6.42}
\elll_2\leq|\pp_{\bar\yy}-\Phi_*\pp_{\bar\yy}|_{var},
\ee
where $\Phi$ is the transformation constructed in Subsection 4.2, and we used the fact that for $\tilde \tau=\iin$ we have $y^v\equiv \xi_v$ and $y^{u'}\equiv \xi_{u'}$
Further, in view \eqref{6.43} we have
\be\label{6.44}
\elll_1\leq \pp_{\bar\yy}(\tau^u\wedge\tau^{u'}<\iin)+\pp_{\bar\yy}(\tau^{u'}<\iin)+|\pp_{\bar\yy}-\Phi_*\pp_{\bar\yy}|_{var}.
\ee
Combining inequalities \eqref{6.45}-\eqref{6.44}, we get
\be\label{6.46}
\pp_\yy(\qqq''_k)\leq 2\,\e_\yy[\ch_{\sigma\geq kT}\cdot (\pp_{\bar\yy}(\tau<\iin)+|\pp_{\bar\yy}-\Phi_*\pp_{\bar\yy}|_{var})].
\ee
Let us note that for any $\om\in\{\sigma\geq kT\}$ we have 
\be\label{6. In-6.4,1}
|\ees_{\tilde u}(kT)|\vee|\ees_{\tilde u'}(kT)|\leq |\ees_{\tilde u}(0)|\vee|\ees_{\tilde u'}(0)|+(L+M)kT+r.
\ee
Moreover, it follows from Proposition \ref{4.13} (see the derivation of \eqref{7.8}) that for any $\es>0$ there is $N$ depending only on $\es, \al, L$ and $M$, such that for all $kT\leq t\leq\tau\wedge\tau^{\tilde v}$, on the set $\sigma\geq kT$, we have 
\begin{align} 
|\xi_{\tilde v}(t)-\xi_{\tilde u}(t)|_\h^2&\leq\exp(-\al(t-kT)+\theta)|\xi_{\tilde u}(kT)-\xi_{\tilde u'}(kT)|_\h^2\notag\\
&\leq\exp(-\al(t-kT)/2+\theta)|\xi_{\tilde u}(kT)-\xi_{\tilde u'}(kT)|_\h^2\label{6. In-6.5},
\end{align}
where we set
$$
\theta=\es\cdot(|\ees_{\tilde u}(kT)|\vee|\ees_{\tilde u'}(kT)|+r).
$$
By the same argument as in the derivation of \eqref{4.19}, we have 
\begin{align}\label{6. In-6.6}
\e_\yy(\ch_{\sigma\geq kT}\cdot|\pp_{\bar\yy}-\Phi_*\pp_{\bar\yy}|_{var})&\equiv \e_\yy(\ch_{\sigma\geq kT}\cdot|\pp_{\yy(kT)}-\Phi_*\pp_{\yy(kT)}|_{var})\notag\\
&\leq\f{1}{2}\left(\left(\e_\yy\exp\left[6\max_{1\leq j\leq N}b_j^{-1}\kkk\right]\ch_{\sigma\geq kT}\right)^{\f{1}{2}}-1\right)^{\f{1}{2}},
\end{align}
where 
\begin{align}
\kkk&=C_1\int_0^\infty \{\exp({-\al(t-kT)/2+\theta})|\xi_{\tilde u}(kT)-\xi_{\tilde u'}(kT)|_\h^2\notag\\
&\q\cdot(1+|\ees_{\tilde u}(kT)|\vee|\ees_{\tilde u'}(kT)|+(L+M)t+r)^2\}\,dt\notag\\
&\leq C_2\int_0^\infty \{\exp({-\al(t-kT)/2+\theta})e^{-\al kT}|y-y'|_\h^2\notag\\
&\q\cdot(1+|\ees_{\tilde u}(0)|\vee|\ees_{\tilde u'}(0)|+(L+M)kT+(L+M)t+r)^2\}\,dt\label{6. In-6.7},
\end{align}
and we used inequalities \eqref{6. In-6.4,1}-\eqref{6. In-6.5} combined with the fact that the mean value is taken along the characteristic of the set $\{\sigma\geq kT\}$.
Now let us fix $\es=\es(\al,L,M)>0$ such that $\al/4\geq\es\cdot(L+M)$, and let $C(\al,L, M)>0$ be so large that for any $k\geq 0$ and any $T\geq 1$
$$
\exp(-\al kT/4)(1+(L+M)kT)^2\leq C(\al,L,M).
$$
Combining this inequality with \eqref{6. In-6.7}, we obtain 
\begin{align}
\kkk&\leq C_3\cdot C(\al,L,M)e^{-\f{\al}{4} kT}|y-y'|_\h^2\int_0^\infty e^{-\al t+\es r}(1+(L+M)t+r)^2\,dt\notag\\
&=C(\al,r,L,M)e^{-\f{\al}{4} kT}|y-y'|_\h^2\label{4.27}.
\end{align}
Now recall that $N$ depends only on $\es, L$ and $M$, and $\es$ depends on $\al, L$ and $M$. It follows that $N$ depends only on $\al,L$ and $M$. Let us choose $d=d(\al,r,L,M)>0$ so small that 
\be\label{4.26}
6\max_{1\leq j\leq N}b_j^{-1}C(\al,r,L,M)d\leq 1.
\ee
Then, by inequalities \eqref{6. In-6.6} and \eqref{4.27}, we have
\be\label{6.47}
\e_\yy(\ch_{\sigma\geq kT}\cdot|\pp_{\bar\yy}-\Phi_*\pp_{\bar\yy}|_{var})\leq e^{-\f{\al}{8}kT}|y-y'|_\h\leq e^{-2(k+2)},
\ee
for $T\geq 16\al^{-1}$ and $d\leq e^{-4}$.
Further, by the Markov property and inequality \eqref{6.48} we have
\begin{align*}
e^{-2(k+2)}&\geq \pp_\yy\{kT\leq \tau<\infty\}=\e_\yy[\e_\yy(\ch_{kT\leq\tau<\iin}|\FFF_{kT})]
=\e_\yy[\pp_{\yy(kT)}(\tau<\iin)]\notag\\
&\geq \e_\yy[\ch_{\sigma\geq kT}\cdot\pp_{\yy(kT)}(\tau<\iin)]\equiv \e_\yy[\ch_{\sigma\geq kT}\cdot\pp_{\bar\yy}(\tau<\iin)].
\end{align*}
Combining this inequality with \eqref{6.46} and \eqref{6.47} we obtain
$$
\pp_\yy(\qqq''_k)\leq 4 e^{-2(k+2)}\leq e^{-2(k+1)}.
$$
\ep

\medskip
\nt
Now we are ready to establish \eqref{6.3}-\eqref{6.5}. We have
$$
\pp_{\yy}\{\sigma=\infty\}\geq1-\sum_{k=0}^\infty\pp_\yy\{kT\leq\sigma\leq(k+1)T\}\geq\f{1}{2},
$$
where used Lemma \ref{6.27} to show that
\begin{align*}
&\pp_\yy\{kT\leq\sigma\leq(k+1)T\}=\pp_\yy\{kT\leq\tau\leq (k+1)T,\tau\leq\vo\}\\
&\q+\pp_\yy\{kT\leq\vo\leq (k+1)T,\vo<\tau\}=\pp_\yy(\qqq_k')+\pp_\yy(\qqq_k'')\leq e^{-2(k+1)}.
\end{align*}
By the same argument,
\begin{align*}
&\e_\yy[\ch_{\{\sigma<\infty\}} e^{\De\sigma}]=\e_\yy[\ch_{\{\sigma<\infty,\tau\leq\vo\}} e^{\De\sigma}]+\e_\yy[\ch_{\{\sigma<\infty,\vo<\tau\}} e^{\De\sigma}]\\
&\leq \e_\yy[\ch_{\{\tau<\infty,\tau\leq\vo\}} e^{\De\tau}]+\e_\yy[\ch_{\{\vo<\infty,\vo<\tau\}} e^{\De\vo}]\leq 2\sum_{k=0}^\infty e^{-2(k+1)}e^{\De k(T+1)}\leq 2,
\end{align*}
for $\De<(1+T)^{-1}$.
So, inequalities \eqref{6.3} and \eqref{6.4} are established. To prove \eqref{6.5}, note that in view of \eqref{6 In-6,1,1}, for $\sigma<\infty$ we have 
\be\label{6.64}
|\ees_{\tilde u}(\sigma)|\leq |\ees_{\tilde u}(0)|+(L+M)\sigma+r.
\ee
Combining this inequality with \eqref{7.9} we obtain
$$
|S_{\sigma}|_{\h}^8\leq   (2|\ees_{\tilde u}(\sigma)|+4C)^4\leq(2(|\ees_{\tilde u}(0)|+(L+M)\sigma+r+2C))^4 \leq C(r,L,M)(1+\sigma^4).
$$
It is clear that the above inequality is satisfied also for $S$ replaced by $S'$, so that
$$
|\SSS_{\sigma}|_{\hh}^{8}\leq 2C(r,L,M)(1+\sigma^4).
$$
Multiplying this inequality by $\ch_{\{\sigma<\infty\}}$, taking the $\e_\yy$-mean value, and using inequality \eqref{6.4}, we arrive at \eqref{6.5}. The proof of Theorem \ref{2.14} (and with it of Theorem \ref{2.15})  is complete.

\subsection{Relaxed non-degeneracy condition}
We finish this section with the following result that allows to relax the non-degeneracy condition imposed on the force.
\begin{theorem}\label{6.10}
There exists $N$ depending only on $\gamma, f, \|h\|$ and $\BBB$ such that the conclusion of Theorem \ref{2.15} remains true for any random force of the form \eqref{1.9}, whose first $N$ coefficients $b_j$ are not zero.
\end{theorem}
Let us fix an integer $N_1$ such that inequality \eqref{4.27} holds for any $N\geq N_1$ and let $d=d(N_1)$ be so small that we have \eqref{4.26}, where $N$ should be replaced by $N_1$. 
Theorem \ref{6.10} will be established, if we show that there is an $N=N(d)\geq N_1$ such
that inequality \eqref{6.13} holds, provided $b_j\neq 0$ for $j=1,\ldots, N$. 
Let the constants $T$ and $\es$, together with the process $\tilde\zeta^y$ be the same as in the proof of Proposition \ref{6.12}. 
We need the following lemma, which is established in the appendix.
\begin{lemma}\label{6.25}
There exists $N\geq N_1$ such that 
\begin{align}\label{6.37}
\sup_{y\in B_R}|(I-P_N)\tilde\zeta^y|_{C(0,T;L^2)}\leq\es/4.
\end{align}
\end{lemma}
We claim that Theorem \ref{6.10} holds with this $N$. Indeed, let us suppose that inequality \eqref{6.13} does not hold, and we have \eqref{6.58}. Let $\tilde\zeta_j$ and $\tilde\zeta$ be the processes constructed in Proposition \ref{6.12}. Denote $\ccc=C(0,T;L^2)$. Then
\begin{align*}
\pp(|\zeta-P_N\tilde\zeta|_{\ccc}\leq \es/4)&=\pp(|\zeta-P_N\tilde\zeta_j+P_N\tilde\zeta_j-P_N\tilde\zeta|_{\ccc}\leq\es/4)\\
&\leq \pp(|\zeta-P_N\tilde\zeta_j|_{\ccc}\leq |P_N\tilde\zeta_j-P_N\tilde\zeta|_{\ccc}+\es/4)\\
&\leq \pp(|\zeta-\tilde\zeta_j+(I-P_N)\tilde\zeta_j|_{\ccc}\leq |\tilde\zeta_j-\tilde\zeta|_{\ccc}+\es/4)\\
&\leq \pp(|\zeta-\tilde\zeta_j|_{\ccc}\leq |(I-P_N)\tilde\zeta_j|_{\ccc}+3\es/4)\leq \pp(|\zeta-\tilde\zeta_j|_{\ccc}\leq \es).
\end{align*}
Letting $j$ go to infinity, and using inequality \eqref{6.22} we obtain
$$
\pp(|\zeta-P_N\tilde\zeta|_{\ccc}\leq \es/4)=0,
$$
which is impossible, since the support of $\zeta$ restricted to $[0,T]$ contains $C(0,T;P_N L^2)$. The proof Theorem \ref{6.10} is complete.

\bigskip
  \section{Appendix}
  
\medskip
\subsection{Proof of \eqref{6.40}} 
Let us consider the continuous map $\gi$ from $C(0,T;H^1_0(D))$ to $C(0,T;\h)$ defined by $\gi(\ph)=\tilde y$, where $\tilde y$ is the flow of equation
$$
\p_t^2 z+\gamma \p_t z-\de z+f(z)-P_N f(z)=h(x)+\p_t \ph,\q [z(0),\dt z(0)]=y.
$$
Then for any $\Gamma\in\bbb(C(0,T;\h))$, we have
\begin{align*}
\pp\{\xi_{\tilde u}(t)\in\Gamma\}&=\pp\{\gi({\int_0^t \psi(s)\,ds})\in\Gamma\}=\pp\{\int_0^t \psi(s)\,ds\in\gi^{-1}(\Gamma)\}\\
&=\pp\{\zeta(t)-\int_0^t P_N f(u)\,ds\in\gi^{-1}(\Gamma)\}\\
&=\pp\{\gi({\zeta(t)-\int_0^t P_N f(u)\,ds})\in\Gamma\}=\pp\{\xi_{u}(t)\in\Gamma\}.
\end{align*}

\medskip
\subsection{Proof of lemma \ref{6.20}} 
Let $f$ be a function that satisfies the growth restriction \eqref{5.17} with $\rho<2$. We claim that inequality \eqref{6.39} holds with $s=(2-\rho)/(2(\rho+1))$. Indeed, by the H\"older and Sobolev inequalities, we have 
\begin{align*}
|f(u)-f(v)|^2_{L^2}&=\int |f(u)-f(v)|^2\leq C\int (|u|^{2\rho}+|v|^{2\rho}+1)|u-v|^2\\
&\leq C||u|^{2\rho}+|v|^{2\rho}+1|_{L^{3\rho/(1-s)}}|u-v|^2_{L^{6/(1+2s)}}\\
&\leq C'(|u|_{H^{1-s}}^{2\rho}+|v|_{H^{1-s}}^{2\rho}+1)|u-v|^2_{H^{1-s}}.
\end{align*}

\medskip 
\subsection{Proof of lemma \ref{6.25}}
First let us fix $N_2\geq N_1$ such that for any $N\geq N_2$
\be\label{6.38} 
|(I-P_N)th|_{C(0,T;L^2)}= T\|(I-P_N)h\|\leq\es/8.
\ee
Notice that
\be\label{6.26}
|(I-P_N)\int_0^t f(\bar u)\,ds |_{C(0,T;L^2)}\leq |I-P_N|_{L(H^{1,p}\to L^2)}T\sup_{0\leq t\leq T}|f(\bar u)|_{H^{1,p}},
\ee
where $p=6/(3+\rho)>6/5$. 
Using the H\"older and Sobolev inequalities, we obtain
\be\label{6.36}
|f(\bar u)|_{H^{1,p}}\leq C_1(|f(\bar u)|_{L^2}+\sum_{j=1}^3|\p_j \bar u f'(\bar u)|_{L^p})\leq C_2(1+\|\bar u\|_1^3)\leq C(R).
\ee
As was already mentioned, the space $H^{1,p}(D)$ is compactly embedded in $L^2(D)$ for $p>6/5$, so that the sequence  $|I-P_N|_{L(H^{1,p}\to L^2)}$ goes to zero as $N$ goes to infinity. Combining this with inequalities \eqref{6.38} -\eqref{6.36}, we arrive at \eqref{6.37}.

\addcontentsline{toc}{section}{Bibliography}
\bibliographystyle{plain}

\begin{thebibliography}{}

\end{thebibliography}


\begin{thebibliography}{10}

\bibitem{BV1992}
A.~V. Babin and M.~I. Vishik.
\newblock {\em Attractors of {E}volution {E}quations}.
\newblock North-Holland Publishing, Amsterdam, 1992.

\bibitem{BCK2014}
Y~Bakhtin, E.~Cator, and K~Khanin.
\newblock Space-time stationary solutions for the {B}urgers equation.
\newblock {\em J. Amer. Math. Soc.}, 27(1):193--238, 2014.

\bibitem{BD-2002}
V.~Barbu and G.~{Da Prato}.
\newblock The stochastic nonlinear damped wave equation.
\newblock {\em Appl. Math. Optim.}, 46(2-3):125--141, 2002.

\bibitem{BKL-2002}
J.~Bricmont, A.~Kupiainen, and R.~Lefevere.
\newblock Exponential mixing of the 2{D} stochastic {N}avier--{S}tokes
  dynamics.
\newblock {\em Comm. Math. Phys.}, 230(1):87--132, 2002.

\bibitem{cazenave1998introduction}
T.~Cazenave and A.~Haraux.
\newblock {\em An introduction to semilinear evolution equations}, volume~13.
\newblock Oxford University Press, 1998.

\bibitem{CV2002}
V.~V. Chepyzhov and M.~I. Vishik.
\newblock {\em {A}ttractors for {E}quations of {M}athematical {P}hysics},
  volume~49 of {\em AMS Coll. Publ.}
\newblock AMS, Providence, 2002.

\bibitem{DZ1992}
G.~{Da Prato} and J.~Zabczyk.
\newblock {\em Stochastic {E}quations in {I}nfinite {D}imensions}.
\newblock Cambridge University Press, Cambridge, 1992.

\bibitem{DZ1996}
G.~{Da Prato} and J.~Zabczyk.
\newblock {\em Ergodicity for {I}nfinite {D}imensional {S}ystems}.
\newblock Cambridge University Press, Cambridge, 1996.

\bibitem{debussche2013ergodicity}
A.~Debussche.
\newblock Ergodicity results for the stochastic {N}avier--{S}tokes equations:
  an introduction.
\newblock In {\em Topics in Mathematical Fluid Mechanics}, pages 23--108.
  Springer, 2013.

\bibitem{DO-2005}
A.~Debussche and C.~Odasso.
\newblock Ergodicity for a weakly damped stochastic non-linear {S}chr\"odinger
  equation.
\newblock {\em J. Evol. Equ.}, 5(3):317--356, 2005.

\bibitem{debussche2013invariant}
A.~Debussche and J.~Vovelle.
\newblock Invariant measure of scalar first-order conservation laws with
  stochastic forcing.
\newblock {\em arXiv preprint arXiv:1310.3779}, 2013.

\bibitem{DirSoug2005}
N.~Dirr and P.~Souganidis.
\newblock Large-time behavior for viscous and nonviscous {H}amilton-{J}acobi
  equations forced by additive noise.
\newblock {\em SIAM J. Math. Anal.}, 37(3):777--796 (electronic), 2005.

\bibitem{EMS-2001}
W.~E, J.~C. Mattingly, and Ya. Sinai.
\newblock Gibbsian dynamics and ergodicity for the stochastically forced
  {N}avier--{S}tokes equation.
\newblock {\em Comm. Math. Phys.}, 224(1):83--106, 2001.

\bibitem{EWKMS2000}
Weinan E, K.~Khanin, A.~Mazel, and Ya. Sinai.
\newblock Invariant measures for {B}urgers equation with stochastic forcing.
\newblock {\em Ann. of Math. (2)}, 151(3):877--960, 2000.

\bibitem{FM-1995}
F.~Flandoli and B.~Maslowski.
\newblock Ergodicity of the 2{D} {N}avier--{S}tokes equation under random
  perturbations.
\newblock {\em Comm. Math. Phys.}, 172(1):119--141, 1995.

\bibitem{MR1641664}
T.~Girya and I.~Chueshov.
\newblock Inertial manifolds and stationary measures for stochastically
  perturbed dissipative dynamical systems.
\newblock {\em Mat. Sb.}, 186, 1995.

\bibitem{HM-2008}
M.~Hairer and J.~C. Mattingly.
\newblock Spectral gaps in {W}asserstein distances and the 2{D} stochastic
  {N}avier--{S}tokes equations.
\newblock {\em Ann. Probab.}, 36(6):2050--2091, 2008.

\bibitem{Har85}
A.~Haraux.
\newblock Two remarks on hyperbolic dissipative problems.
\newblock In {\em Nonlinear partial differential equations and their
  applications. {C}oll\`ege de {F}rance seminar, {V}ol.\ {VII} ({P}aris,
  1983--1984)}, volume 122 of {\em Res. Notes in Math.}, pages 6, 161--179.
  Pitman, Boston, MA, 1985.

\bibitem{IturK2003}
R.~Iturriaga and K.~Khanin.
\newblock Burgers turbulence and random {L}agrangian systems.
\newblock {\em Comm. Math. Phys.}, 232(3):377--428, 2003.

\bibitem{kuksin2013stochastic}
S.~Kuksin and V.~Nersesyan.
\newblock Stochastic {CGL} equations without linear dispersion in any space
  dimension.
\newblock {\em Stochastic Partial Differential Equations: Analysis and
  Computations}, 1(3):389--423, 2013.

\bibitem{KS-cmp2000}
S.~Kuksin and A.~Shirikyan.
\newblock Stochastic dissipative {PDE}s and {G}ibbs measures.
\newblock {\em Comm. Math. Phys.}, 213(2):291--330, 2000.

\bibitem{KS-book}
S.~Kuksin and A.~Shirikyan.
\newblock {\em Mathematics of {T}wo-{D}imensional {T}urbulence}.
\newblock Cambridge University Press, Cambridge, 2012.

\bibitem{Lions1969}
J.-L. Lions.
\newblock {\em Quelques m\'ethodes de r\'esolution des probl\`emes aux limites
  non lin\'eaires}.
\newblock Dunod; Gauthier-Villars, Paris, 1969.

\bibitem{MR1245306}
C.~Mueller.
\newblock Coupling and invariant measures for the heat equation with noise.
\newblock {\em Ann. Probab.}, 21(4):2189--2199, 1993.

\bibitem{odasso-2008}
C.~Odasso.
\newblock Exponential mixing for stochastic {PDE}s: the non-additive case.
\newblock {\em Probab. Theory Related Fields}, 140(1-2):41--82, 2008.

\bibitem{shirikyan-ptrf2006}
A.~Shirikyan.
\newblock Law of large numbers and central limit theorem for randomly forced
  {PDE}'s.
\newblock {\em Probab. Theory Related Fields}, 134(2):215--247, 2006.

\bibitem{shirikyan-bf2008}
A.~Shirikyan.
\newblock Exponential mixing for randomly forced partial differential
  equations: method of coupling.
\newblock In {\em Instability in models connected with fluid flows. {II}},
  volume~7 of {\em Int. Math. Ser. (N. Y.)}, pages 155--188. Springer, New
  York, 2008.

\bibitem{VKF-1979}
M.~I. Vishik, A.~I. Komech, and A.~V. Fursikov.
\newblock Some mathematical problems of statistical hydromechanics.
\newblock {\em Uspekhi Mat. Nauk}, 34(5(209)):135--210, 1979.

\bibitem{zelik2004}
S.~Zelik.
\newblock Asymptotic regularity of solutions of a nonautonomous damped wave
  equation with a critical growth exponent.
\newblock {\em Commun. Pure Appl. Anal}, 3:921--934, 2004.

\end{thebibliography}
\def\cprime{$'$} \def\cprime{$'$}
  \def\polhk#1{\setbox0=\hbox{#1}{\ooalign{\hidewidth
  \lower1.5ex\hbox{`}\hidewidth\crcr\unhbox0}}}
  \def\polhk#1{\setbox0=\hbox{#1}{\ooalign{\hidewidth
  \lower1.5ex\hbox{`}\hidewidth\crcr\unhbox0}}}
  \def\polhk#1{\setbox0=\hbox{#1}{\ooalign{\hidewidth
  \lower1.5ex\hbox{`}\hidewidth\crcr\unhbox0}}} \def\cprime{$'$}
  \def\polhk#1{\setbox0=\hbox{#1}{\ooalign{\hidewidth
  \lower1.5ex\hbox{`}\hidewidth\crcr\unhbox0}}} \def\cprime{$'$}
  \def\cprime{$'$} \def\cprime{$'$} \def\cprime{$'$}

\end{document}